\documentstyle[12pt,sr-stylefile]{article}

\begin{document}

\def\f{\footnote}
\def\Q{\begin{quote} \baselineskip=12pt}
\def\q{\end{quote} \baselineskip=24pt}
\renewcommand{\thefootnote}{}
\def\sgn{{\rm sgn}}
\def\id{{\rm Id}}
\def\even{{\rm even}}
\def\Pf{{\rm Pf}}
\def\dim{{\rm dim}}
\def\dvol{{\rm dvol}}
\def\MQ{{\rm MQ}}
\newcommand{\bbC}{C}
\newcommand{\bbR}{{\bf R}}
\def\evev{(\exp_{\bar x}v,\exp_{\bar x} (-v))} 
\def\dxy{d_{x,y}}
\def\cutx{{\cal C}_x}
\def\expx{\exp_{\bar x}}
\def\adel{A_\delta}
\def\bdel{B_\delta}
\def\itfmq{({\rm Id},tf)^*{\rm MQ}_\Delta}
\def\euler{{\rm Pf}(\Omega)}
\def\ade{A_{\delta,\epsilon}}
\def\bde{B_{\delta,\epsilon}}
\def\cutf{{\cal C}(f)}
\def\degf{{\rm deg}(f)}
\def\cf{{\cal C}^f}
\def\vfp{V_f^\prime}
\def\tmup{TM^\uparrow}
\def\itfup{({\rm Id},tf)^*{\rm MQ}_{TM^\uparrow}}
\def\mqup{{\rm MQ}_{TM^\uparrow}}
\newcommand{\ALpha}{\frac{e^{-\rho^2(\frac{d(x,f(x))}{\sqrt{2}})}}
{(2\pi)^{n/2}}}
\newcommand{\rats}{\frac{\epsilon (I,I^\prime)}{|I|!\ |I^\prime|!}}
\def\deg{{\rm deg}}
\def\vol{{\rm vol}}

\Large
\centerline{\bf Mathai-Quillen Forms and Lefschetz Theory}
\bigskip\bigskip
\large
\centerline{{\bf Mihail Frumosu}${}^\dag$
}
\normalsize 
\medskip
\centerline{\parbox{5in}{
 Department of Mathematics, Boston University,
Boston, MA 02215}}
\centerline{ matf@math.bu.edu}\footnote{\dag \ Current address:
Metropolitan College, Boston University, Boston,
MA 02215}

\bigskip
\large
\centerline{{\bf Steven Rosenberg}${}^\ddag$\footnote{\ddag\ The
 second
author's research was partially
supported by the NSF.} }
\normalsize
\medskip
\centerline{\parbox{5in}{
 Department of Mathematics, Boston University,
Boston, MA 02215}}
\centerline{ sr@math.bu.edu}

\bigskip
\bigskip

\section{\large\rm\bf Introduction}

In \cite{MQ}, Mathai and Quillen introduced a
geometric representative, the Mathai-Quillen form, for the Thom class
of an oriented Riemannian bundle over an oriented Riemannian
manifold.  Using a one-parameter family of
pullbacks of this form, Mathai and Quillen gave a new proof of both
the Hopf index formula and the Chern-Gauss-Bonnet formula for the
Euler characteristic.

The Euler characteristic is the Lefschetz number of the identity map
of a closed manifold $M$, and various topological
expressions for the Euler characteristic 
(alternating sum of Betti numbers,
self-intersection number of the diagonal in $M\times M$, Hopf index formula)
have counterparts in Lefschetz theory for functions 
$f:M\to M$ (supertrace of $f$ on
cohomology, intersection number of $f$'s graph with the diagonal, Lefschetz
fixed point formula). However, a Lefschetz counterpart of the
integral 
geometric expression for the Euler characteristic (the
Chern-Gauss-Bonnet
theorem) has not been known previously.  With this motivation, we
investigate 
geometric aspects of
Lefschetz theory 
using Mathai-Quillen forms, and in particular study the analogous
one-parameter family of pullback forms.

In \S2.2, we give via Poincar\'e duality
an elementary integral formula for the Lefschetz
number in terms of the map $f$ and the Mathai-Quillen form of the
normal bundle of the diagonal in $M\times M$ (Theorem
\ref{theoremone}).  This formula
specializes to the Chern-Gauss-Bonnet formula when $f = \id$.  
In contrast to the Chern-Gauss-Bonnet integrand, the local expression
of the
integrand is fairly complicated even for flat manifolds, due to the
action of $f$ on the Mathai-Quillen form.  In \S2.3, we compute the
integrand in the flat case (Theorem \ref{theoremtwo})
and work a computation on $S^1.$
We give the integrand for general metrics in \S2.4 (Theorem
\ref{localtheorem}).  Here the result
is less explicit, as the integrand depends on solutions of Jacobi
fields along geodesics joining $x$ and $f(x).$  These Jacobi fields
appear because we must exponentiate the Mathai-Quillen form of the
normal bundle  to a form on a tubular neighborhood of
the diagonal,
and the Jacobi fields measure the deviation of the
exponential map from the trivial flat case.  We are able to make the
integrand
explicit for constant curvature metrics.

In \S\S3-4, we investigate the Lefschetz analog of the one-parameter
family of pullbacks of Mathai-Quillen
forms.  In brief, the $t\to\infty$ limit of the pullbacks
measures the
fixed point set of the map, while the $t\to 0$ limit measures the set
of points mapped far from themselves.

More precisely, as $t\to\infty$, we recover in \S3 the
Lefschetz fixed point formula (Theorem \ref{theoremapp}), 
and in fact we give a new proof of
the more general formula for submanifolds of fixed points.  This
argument is done both at the topological level using
Thom classes, and at the
geometric level using Mathai-Quillen forms.  The proof can be thought
of as a simplified version 
of the heat equation proof of the Lefschetz formula in
\cite{G}.

In \S4, we consider the $t\to 0$ limit.  This
limit is trivial in the case considered in \cite{MQ}, but is
discontinuous in our case.  In fact, the discontinuities occur at the
intersection of the graph of $f$ with the boundary of the tubular
neighborhood of the diagonal, which in general is quite complicated.  

To extract the
maximum geometric information, we take the largest possible tubular
neighborhood, with boundary  the union of the cut loci of
the points of $M$.   This choice is motivated by the observation that
if $f(x)$ is never in $\cutx$, the cut locus of $x$, then $L(f)  =
\chi(M).$
More precisely, there is a singular current supported on $\cutf =
\{x:f(x)
\in\cutx\}$, whose singular part measures $L(f) - \chi(M)$ (Theorem
\ref{singcur}).  Assuming that $\cutf$ is finite and
imposing a transversality condition,
we find the sharp estimate
$|L(f) - \chi(M)| \leq |\cutf|$ 
 (Theorem \ref{isolated}).  
These assumptions place strong restrictions on the metric on $M$, and 
for diffeomorphisms $f$ with $L(f) \neq \chi(M)$, $|\cutf|$ is
infinite for most metrics (Theorem \ref{yahthree}).

These results give geometric information for Lefschetz theory via
Mathai-Quillen forms.  
In the
appendix, we show how Hodge theory techniques give upper bounds for the
Lefschetz number in terms of the geometry of $M$.

\section{\large\rm\bf The basic formula and its local
expression}\label{firstsection} 

Let $f:M\to M$ be a smooth map of a closed oriented
Riemannian manifold $M$.  After a review of Poincar\'e duality in
\S2.1, 
we give in \S2.2 an integral formula (Theorem \ref{theoremone}) 
for the Lefschetz number of
$f$ which reduces to the Chern-Gauss-Bonnet theorem when $f = {\rm
Id}$.  The Mathai-Quillen formalism is also easily
extended to odd rank bundles.
The local expression for the integrand for flat manifolds is
computed in \S2.3 (Theorem
\ref{theoremtwo}), and
the integrand for
arbitrary metrics is computed in \S2.4 (Theorem \ref{localtheorem}).  
This last formula is
then specialized to constant curvature metrics.

\subsection{\large\rm\bf  Topological preliminaries}

The Lefschetz number of $f$ is
$$L(f)=\sum_q(-1)^q\ {\rm tr}\ f^{q},$$
where $f^{q}$ denotes the induced map on the real cohomology
group $H^q(M)$. 
This has a well known Poincar\'e duality formulation,
\cite[Ex.~11.26]{BT}, which we give for completeness.

\begin{lemma} \label{firstlemma}
Let $f:M\rightarrow M$ be a smooth map of
a closed oriented manifold. Then
$$L(f)=\int_\Delta\eta_\Gamma,$$
where $\eta_\Gamma$ is the Poincar\'e dual of the graph $\Gamma$ of
$f$
in $M\times M$ and $\Delta$ is the diagonal of $M$.
\end{lemma}

Here we do not distinguish between the cohomology class $\eta_\Gamma$
and a representative form.  Before the proof of the lemma, we collect
the basic results about Poincar\'e duality and Thom classes.  Recall
that the Poincar\'e dual $\eta_N$ 
of an oriented $k$-submanifold $N$ of a closed oriented manifold $X$ is
the real cohomology class
defined by (or characterized by, depending on one's definition)
\begin{equation}\label{pduality} \int_N \omega = \int_X
\omega\wedge\eta_N,\end{equation} 
for all closed $k$-forms $\omega$ on $X$ \cite[(5.13)]{BT}.  
\begin{theorem}  \label{bttheorem}
(i) Let $N'$ be another closed oriented submanifold
of $X$ with transverse intersection with $N$. Then
\[ \eta_{N\cap N'} = \eta_N\wedge\eta_{N'}.\]

(ii)  A closed form $U\in H_c^k(E)$, the compactly supported
cohomology of an oriented  rank $k$ bundle $E$ over $X$, represents
the Thom
class iff the integral of $U$ over each fiber of $E$ is one.

(iii)  Identify the total space of $\nu_N^X$, the
normal bundle of $N$ in $X$, with a tubular
neighborhood of $N$ in $X$, so that the Thom class of the normal
bundle
can be considered as a cohomology class on $X$.  Then
the Poincar\'e dual of $N$ is the same as the
Thom class of the normal bundle of $N$ in $X$. \end{theorem}

The proofs of (i)-(iii) are in (6.31), Prop. 6.18 and Prop.~6.24 of
\cite{BT}, respectively.  It is pointed out in \cite{MQ} that
the cohomology with compact support in (ii) may be replaced with the
cohomology of forms with $C^1$ exponential decay in the fibers.

\bigskip

\noindent {\bf Remark:}
The main technical work in this paper is in making the identification
in (iii) explicit.  To consider a closed form $U$ on $\nu_N^X$ as a
closed form on $X$, 
we first use a diffeomorphism
$\alpha$
from
the $\epsilon$-ball $B_\epsilon (0)\subset {\bf R}^n$
to ${\bf R}^n$
to pull $U$ back to the form $\alpha^*U$ on the
$\epsilon$-neighborhood of the zero section in $\nu_N^X$.
For $\epsilon$ small enough, the exponential map $\exp: \nu_N^X\to
X$ is a diffeomorphism onto a tubular neighborhood
of $N$, and so it is really
$(\exp^{-1})^*\alpha^*U$ which is a form supported on the tubular
neighborhood.

\bigskip

\noindent {\sc Proof of Lemma 2.1:} 
Let $\lbrace\omega_i\rbrace$ be a basis for  
$H^*(M)$ and $\lbrace\tau_j\rbrace$ the dual basis under 
Poincar\'e duality, i.e.~
$\int_M\omega_i\wedge\tau_j=\delta_{ij}.$
Let $\pi,\ \rho$ be the projections of $M\times M$ onto the first and
second factors.
$H^*(M\times M)$ has as basis
$\{\pi^*\omega_i\wedge\rho^*\tau_j\}$, so
$\eta_\Gamma=\sum_{i,j}c_{ij}\pi^*\omega_i\wedge\rho^*\tau_j$
for some $c_{ij}\in {\bf R}$. We now determine the $c_{ij}.$

\begin{lemma} \label{secondlemma}
$\displaystyle{\eta_\Gamma=\sum_{i,j} 
(-1)^{(\deg\ \omega_i)(\deg\ \omega_j)}\alpha_{ji}\pi^*\omega_i
\wedge\rho^*\tau_j}$,
with $\alpha_{ij}$ defined by
$f^*\omega_i=\alpha_{ij}\omega_j.$ \end{lemma}

\noindent {\sc Proof}: We compute 
$\int_\Gamma \pi^*\tau_k\wedge\rho^*\omega_l$ in two ways. 
Using the graph map
$i:M\longrightarrow \Gamma\subset M\times M\ , i(x)=(x, f(x))$,
we obtain
\begin{eqnarray*}\int_\Gamma \pi^*\tau_k\wedge\rho^*\omega_l
&=&\int_M i^*\pi^*\tau_k\wedge i^*\rho^*\omega_l
=\int_M \tau_k\wedge f^*\omega_l\cr
&=&\int_M \alpha_{lj}\tau_k\wedge\omega_j
=\alpha_{lj}(-1)^{(\deg\ \omega_j)(\deg\ \tau_k)}\delta_{kj}\cr
&=&\alpha_{lk}(-1)^{(\deg\ \omega_k)(\deg\ \tau_k)}.\cr\end{eqnarray*}
By (\ref{pduality}),
\begin{eqnarray*}\int_\Gamma \pi^*\tau_k\wedge\rho^*\omega_l
&=&\int_{M\times M} \pi^*\tau_k\wedge
\rho^*\omega_l\wedge\eta_\Gamma\cr
&=&\sum_{i,j}c_{ij}\int_{M\times M} \pi^*\tau_k\wedge
\rho^*\omega_l\wedge 
\pi^*\omega_i\wedge\rho^*\tau_j\cr
&=&\sum_{i,j}c_{ij}(-1)^{(\deg\ \tau_k +\deg\ \omega_l)(\deg\
\omega_i)}
\int_{M\times
M}\pi^*(\omega_i\wedge\tau_k)\wedge\rho^*(\omega_l\wedge\tau_j)\cr
&=&(-1)^{(\deg\ \tau_k +\deg\ \omega_l)(deg\
\omega_l)}c_{kl}.\cr\end{eqnarray*}
Thus
$c_{kl}=\alpha_{lk}(-1)^{(\deg\ \omega_k)(\deg\ \omega_l)}.$
\hfill$\Box$
\medskip

For the proof of Lemma \ref{firstlemma}, we have
\begin{eqnarray*}\int_\Delta\eta_\Gamma 
&=&\sum_{i,j}(-1)^{(\deg\ \omega_i)(\deg\ \omega_j)}
\alpha_{ji}\int_M i^*\pi^*\omega_i\wedge 
i^*\rho^*\tau_j\cr
&=&\sum_{i,j}(-1)^{(\deg\ \omega_i)(\deg\ \omega_j)}
\alpha_{ji}\int_M\omega_i\wedge\tau_j
=\sum_{i}(-1)^{(\deg\ \omega_i)}\alpha_{ii}\cr
&=&\sum_q(-1)^q\ {\rm tr}\ f^q
= L(f).\cr\end{eqnarray*}
\hfill$\Box$
\medskip

If $\Gamma$ is transversal 
to $\Delta$ in $M\times M$, this lemma leads to a quick
proof of the Lefschetz fixed point formula:
$$L(f)=\sum_{p, f(p) = p}\sigma_p,$$
with  $\sigma_p=\sgn\ \det(\id -(df)_p)$ and $(df)_p: T_pM 
\to T_pM$ the derivative map.
Here $\Gamma\cap\Delta$
is a finite set of points. 
Theorem \ref{bttheorem} and Poincar\'e duality give
\begin{eqnarray*}L(f)
&=&\int_\Delta\eta_\Gamma
=\int_{M\times M}\eta_{\Gamma}\wedge\eta_{\Delta}
=\int_{M\times M}\eta_{\Gamma\cap\Delta}\cr
&=&\int_{\Gamma\cap\Delta} 1
=\sum_{p, f(p) =p}\pm 1.\cr\end{eqnarray*}
Thus $L(f)$ is the sum of the orientations $\pm 1$ of the fixed
points $p$ of $f.$ By [GP, p.121], the orientation equals
$\sgn\ \det(\id -(df)_p)$ in our sign convention.
Note also that $L(f)=\int_{\Gamma\cap\Delta} 1$ implies
that $L(f)=I(\Delta ,\Gamma)$, the intersection number of $\Delta$ and
$\Gamma$, so we have three equivalent definitions of the
Lefschetz.

\subsection{\large\rm\bf Mathai-Quillen formalism and the integral
formula}

In \cite{MQ}, Mathai and Quillen obtained a geometric expression for
the Thom class of an oriented even dimensional vector bundle.
Let $E$ be a rank $n=2m$ vector bundle over a
manifold
$M$, where $E$ has an inner product and a compatible connection
$\theta$.
Then a geometric 
representative $\MQ$ of the  Thom class of $E$ is given by
\begin{equation}\label{basicdef}
\MQ=\pi^{-m}e^{-x^2}\sum_{I,|I|\ \even}\epsilon
(I,I')\Pf\bigg({1\over 2}\Omega_I\bigg)(dx+\theta x)^{I'},\end{equation}
where: $x$ is an orthornormal fiber coordinate; $\Omega$ 
is curvature of
the connection $\theta$; $\Omega_I$ is the submatrix of $\Omega$ with
respect to the multi-index $I$ with entries in $\{1,2,... ,n\}$;
$\Pf({1\over
2}\Omega_I)$
is the Pfaffian of ${1\over 2}\Omega_I$; $I'$ denotes the complement
of $I$ in $\{1,2,... ,n\}$; $\epsilon (I,I')$ is the sign of $I,I'$
considered as a shuffle permutation in the exterior algebra:
$$dx^I\wedge dx^{I'}=\epsilon (I,I')dx^1\wedge ...\wedge dx^n;$$
and
$$(dx+\theta x)^{I'}
=(dx^{i_1}+\theta^{i_1}_{j_1}x^{j_1})\wedge
(dx^{i_2}+\theta^{i_2}_{j2}x^{j_2})\wedge ...\wedge
(dx^{i_q}+\theta^{i_q}_{j_q}x^{j_q}),$$
with $I'=\{i_1, i_2,... ,i_q\}$. In the expression $\theta x$,
$\theta$ denotes the connection one-forms of the connection for the
frame $\{ x^i\}$. The ordering of the elements of
$I^\prime$ in $dx^{I^\prime}$ is unimportant due to the $\epsilon
(I,I')$ factor.
 For computations at a point $x\in M$, we will often
assume that $\{x^i\}$ is a synchronous frame centered at $x$, in which
case the connection one-forms $\theta$ vanish at $x$

Unlike the Euler characteristic, the Lefschetz number can be nonzero
for odd dimensional manifolds, so
we need to check that this formalism extends to bundles of odd rank.
Let $n= 2m+1$ and let $E_M$ be an oriented
 rank $n$ vector bundle over a manifold $M$ with a connection
compatible with a metric on $E_M$, and 
let $E_{S^1}$ be the trivial bundle with the trivial connection
over $S^1$. Equip 
$E_M\times E_{S^1}$ over $M\times S^1$
with the product connection.
The Mathai-Quillen representative
$\MQ_{M\times S^1}\in H^{2m+2}(M\times S^1)$ 
of the Thom class of $E = E_M\times E_{S^1}$ is given by
$$\MQ_{M\times S^1}=\pi^{-(m+1)}e^{-x^2}\sum_{I, |I|\ \even}\epsilon
(I,I')\Pf\bigg({1\over 2}\Omega_I\bigg)(dx+\theta x)^{I'}\ ,$$
where $\theta$ is the connection one-form with respect to a product
orthonormal frame $\{x^i\}$ of $E$, and $\Omega=\Omega_{M\times S^1}$
is the 
curvature of this connection over  $M\times S^1$.
The curvature matrix for $M\times S^1$ is:
$$\Omega_{M\times S^1}=\pmatrix{\Omega_M & 0\cr
                                0 & 0\cr},$$
where $\Omega_M$ is the curvature matrix of $E_M$.

Recall that for an even-dimensional $k\times k$
matrix $\omega$, the Pfaffian is a homogeneous polynomial
of degree $k/2$ in the entries of $\omega$ characterized
up to  sign by
$\Pf^2(\omega)= \det(\omega).$
If $\Omega_I$ is a submatrix of $\Omega_{M\times S^1}$, then 
$\Pf\Big({1\over 2}\Omega_I\Big)=0,$
unless $\Omega_I$ is a submatrix of $\Omega_M$ itself.
Thus in the definition of $\MQ_{M\times S^1}$ we may assume that
$n+1\not\in I$, where $n+1$ corresponds to the $dt$ variable and $t$
is the coordinate on $S^1$.
Moreover, we have
\begin{eqnarray*}(dx+\theta x)^{I'}
&=&(dx^{i_1'}+\theta^{i_1^\prime}_{j_1}x^{j_1})\wedge
(dx^{i_2'}+\theta^{i_2^\prime}_{j2}x^{j_2})\wedge ...\wedge
(dt +\theta^{n+1}_{j_q}x^{j_q})\cr
&=&(dx+\theta x)^{I_M'}\wedge dt,\cr\end{eqnarray*}
where $I'=I_{M}'\cup\{ n+1\}=\{i_1',... ,i_{q-1}', n+1\}$ 
with $I_M'$  the 
complement of $I-\{ n+1\}$ in $\{ 1,2, ...,n\}$.
Hence $\MQ_{M\times S^1}$ decomposes as follows:
\begin{eqnarray*}\MQ_{M\times S^1}
&=&\pi^{-(m+{1\over 2})}e^{-x^2}\sum_{I_M,|I_M|\
\even}\epsilon(I_M,I_M')
\Pf\bigg({1\over 2}\Omega_{I_M}\bigg)(dx+\theta
x)^{I_M'}\cr
&&\qquad \wedge(\sqrt{\pi})^{-1}e^{-t^2}dt\cr
&=&U_M\wedge U_{S^1},\cr\end{eqnarray*}
where 
$$U_{M}=\pi^{-{n\over 2}}e^{-x^2}\sum_{I_M,|I_M|\
\even}\epsilon(I_M,I_M')
\Pf\bigg({1\over 2}\Omega_{I_M}\bigg)(dx+\theta x)^{I_M'},\ \ 
U_{S^1}={\pi}^{-1/2}e^{-t^2}dt.$$

We want to show that $U_M$ represent the Thom class of $E_M$. 
By Theorem \ref{bttheorem} (iii),
$U_{S^1}$ 
represents the Thom class of $E_{S^1}$.

\begin{lemma}\label{lemmathree} If $E$ is a vector bundle of
odd rank $n$ over $M$, then
$$U_{M}=\pi^{-{n\over 2}}e^{-x^2}\sum_{I,|I|\ \even}\epsilon 
(I,I')\Pf\bigg({1\over 2}\Omega_I\bigg)(dx+\theta x)^{I'}$$
is in the Thom class of $E_M$. \end{lemma}

\noindent{\sc Proof}: We know that
$U_M\wedge U_{S^1}$ represents the Thom class for
$E_M\times E_{S^1}$. 
Since $dU_{S^1} =0$, we have
$0
=d(U_{M\times S^1})
=dU_M\wedge U_{S^1}.$
$U_{S^1}$ is non-zero, so $dU_M=0$. 
Moreover, in each fiber 
$\int_{E_{S^1}}U_{S^1}=1$, so in each fiber
\[1
=
\int_{E_{M}\times E_{S^1}}U_M\wedge U_{S^1}
=\Big(\int_{E_{M}} U_M\Big)\times\Big(\int_{E_{S^1}}U_{S^1}\Big)
=\int_{E_M} U_M.\]
Hence $U_M$ represents the Thom class of $E_M$.\hfill$\Box$

\medskip


The following elementary result is the basic integral formula for the
Lefschetz number. 
\begin{theorem} \label{theoremone}
Let $f:M\rightarrow M$ be a smooth
map of a closed oriented manifold. Let
$\Delta_\epsilon$ be a tubular neighborhood of the diagonal in
$M\times M$ of width $\epsilon$, and let 
$\MQ_{\Delta_\epsilon}$ be the Mathai-Quillen form
of the normal bundle to the diagonal, considered as a form supported
 in 
$\Delta_\epsilon$.
Then the Lefschetz
number is given by
\begin{equation}\label{eqnzero}
L(f)=(-1)^{\dim\ M}\int_M
 (\id,f)^* \MQ_{\Delta_\epsilon},\end{equation}
where $(\id,f): M\rightarrow\Gamma\subset M\times M$ is
the graph map.
\end{theorem}

\noindent {\sc Proof:} By Lemma 2.1 and Poincar\'e duality, we have
\begin{eqnarray*}L(f) 
&=&\int_\Delta \eta_\Gamma
=  \int_{M\times M} \eta_\Gamma \wedge \eta_\Delta\cr
&=&  (-1)^{\dim\ M} \int_\Gamma\eta_\Delta
=(-1)^{\dim\ M} \int_{(\id,f)(M)} \MQ_{\Delta_\epsilon}\cr
&=& (-1)^{\dim\ M}\int_M (\id,f)^*\MQ_{\Delta_\epsilon},\cr
\end{eqnarray*}
since $(\id,f)$ is an orientation preserving diffeomorphism and hence 
of degree 1. \hfill$\Box$

\medskip
This formula generalizes the Chern-Gauss-Bonnet theorem.
The Euler characteristic of an even dimensional Riemannian manifold
$M$ is given by
$$\chi (M)=L(\id)=\int_M (\id,\id)^* \MQ_{\Delta_\epsilon}=\int_M 0^*
\MQ_{TM},$$
since a neighborhood of the zero section in
$TM$ is isomorphic to a tubular neighborhood of $\Delta$ under
an isomorphism taking the zero section  0 to the graph map $({\rm
Id},\id)$ of the identity.
For the Levi-Civita connection $\theta$, we have
$$\MQ_{TM}=\pi^{-n/2}\sum_{|I|\ {\rm even}}\epsilon (I,I')\Pf({1\over
2}\Omega_I)
(dx+\theta x)^{I'},$$
and
$$0^*\MQ_{TM}=\pi^{-n/2}\Pf({1\over 2}\Omega)$$
since $x=0$ on $M$ implies
$0^*(dx+\theta x)^{I'}=0$ if $I'\neq\emptyset$.
Thus we obtain the Chern-Gauss-Bonnet theorem
$$\chi(M)={1\over (2\pi)^{n/2}}\int_M \Pf(\Omega).$$
Similarly, we find $\chi(M)=0$ if $\dim\ M$ is odd.

It is important to note that the support of the integrand in
(\ref{eqnzero}) is $\{x\in M:(x,f(x))\in\Delta_\epsilon\}.$

\subsection{\large\rm\bf Local expressions for flat manifolds}

We will check Theorem \ref{theoremone} on a simple flat example and
derive an integral formula for the Lefschetz number of a general flat
manifold. 
\bigskip

\noindent{\bf Example:}
Let $M=S^1$ and let $f:S^1\to S^1$ be given by
$f(z)=z^n$.  Then $f$ has degree $n$, so $L(f) = 1-n.$

We now construct the Mathai-Quillen form $\MQ_{\Delta_\epsilon}$
of an $\epsilon$-neighborhood
of the diagonal in $S^1\times S^1$, where we fix
$\epsilon = {\pi\over
2\sqrt{2}}$ for convenience.
Let
$$\MQ_{\nu_\Delta}={1\over\sqrt{\pi}}e^{-x^2}dx$$
be the Mathai-Quillen form of  $TS^1$, which we identify with
$\nu_\Delta^{S^1\times S^1}$, the normal bundle of the diagonal in
$S^1\times S^1.$
Let
$\alpha:\Delta_\epsilon \to
\nu^{S^1\times S^1}_{\Delta_{S^1}}$ be the diffeomorphism
$$\alpha(\theta_1,\theta_2)=\Biggl({\theta_1+\theta_2\over 2},
\rho\bigg({\theta_1-\theta_2\over \sqrt{2}}\bigg)\Biggr),$$
where $(\theta_1,\theta_2)$ are the coordinates on $S^1\times S^1$ and
$\rho :\Big(-{\pi\over 2\sqrt{2}},{\pi\over 2\sqrt{2}}\Big)
\to (-\infty,\infty)$ is an
orientation preserving  diffeomorphism given by a fixed odd
function $\rho$.  We set $\rho(x) = \infty,\ -\infty$ if $x>{\pi\over
2\sqrt{2}},\ x< -{\pi\over 2\sqrt{2}}$, respectively.

The condition
${1\over\sqrt{\pi}}\int_{-\infty}^{\infty}e^{-x^2}dx=1$
implies
$${1\over\sqrt{\pi}}\int_{-{\pi\over 2\sqrt{2}}}^{\pi\over 2\sqrt{2}}
\rho '(\theta)e^{-{\rho^2(\theta)}}d\theta=1, \ \ {\rm and} \ 
\int_0^{{\pi\over 2\sqrt{2}}}
\rho '(\theta)e^{-\rho^2(\theta)}d\theta={\sqrt{\pi}\over 2},$$
since the integrand is even.
Hence we have
\[
\MQ_{\Delta_\epsilon}
=\alpha^*(\MQ_{\nu_\Delta})
={1\over\sqrt{2\pi}}e^{-\rho^2({\theta_2-\theta_1\over\sqrt{2}})}
\rho
'\bigg({\theta_2-\theta_1\over\sqrt{2}}\bigg)(-d\theta_1+d\theta_2)
\]
at $(\theta_1 ,\theta_2)$.

The graph of $f$, drawn on $[0,2\pi]\times [0,2\pi],$ consists of $n$
line segments 
$\theta_2=n\theta_1-2(k-1)\pi,\ k=1,2,\ldots, n.$
Since the upper and lower limits of the tubular neighborhoods are
given by
$\theta_2=\theta_1\pm{\pi\over 2},$ it is easy to check that $\Gamma$
is in the tubular neighborhood iff
$${(4k-5)\pi\over 2(n-1)}\leq\theta_1\leq
  {(4k-3)\pi\over 2(n-1)},$$
for $k=2, ...,n-1$,  or
$$0\leq\theta_1\leq{\pi\over2(n-1)}\ ,\ \ \ \ 
  {(4n-5)\pi\over 2(n-1)}\leq\theta_1\leq 2\pi,$$
for the first and last segment, respectively.  Thus
\begin{eqnarray*}(\id,f)^*\MQ_{\Delta_\epsilon}
&=&{1\over\sqrt{2\pi}}e^{-\rho^2({(n-1)\theta -(k-1)2\pi\over
\sqrt{2}})}
\rho '\bigg({(n-1)\theta -(k-1)2\pi\over\sqrt{2}}\bigg)\cr
&&\qquad \cdot (\id,f)^*
(d\theta_2-d\theta_1)\cr
&=&\bigg({n-1\over\sqrt{2\pi}}\bigg)
e^{-\rho^2({(n-1)\theta -(k-1)2\pi\over \sqrt{2}})}
\rho '\bigg({(n-1)\theta -(k-1)2\pi\over \sqrt{2}}\bigg)d\theta
,\cr\end{eqnarray*}
since $f^*d\theta = nd\theta.$  This gives
\begin{eqnarray*}L(f) 
&=&-\int_{S^1}(\id,f)^*\MQ_{\Delta_\epsilon}\cr
&=&-{1\over\sqrt{\pi}}\int_0^{\pi\over 2(n-1)}
e^{-\rho^2({{(n-1)\theta\over\sqrt{2}}})}
\rho '\bigg({(n-1)\theta\over\sqrt{2}}\bigg)d\theta\cr
&&\qquad -\sum_{k=2}^{n-1}\biggl[
{1\over\sqrt{\pi}}\int_{(4k-5)\pi\over 2(n-1)}^{(4k-3)\pi\over 2(n-1)}
\bigg({n-1\over\sqrt{2}}\bigg)
e^{-\rho^2\big({{(n-1)\theta\over\sqrt{2}}-(k-1)\pi\sqrt{2}}\big)}\\
&&\qquad\qquad
\rho '\bigg({(n-1)\theta -(k-1)2\pi\over \sqrt{2}}\bigg)
d\theta\biggr]\cr
&&\qquad
-{1\over\sqrt{\pi}}\int_{(4n-5)\pi\over 2(n-1)}^{2\pi}
\bigg({n-1\over\sqrt{2}}\bigg)
e^{-\rho^2\big({(n-1)(\theta -2\pi)\over \sqrt{2}}\big)}
\rho '\bigg({(n-1)(\theta -2\pi)\over \sqrt{2}}\bigg)
d\theta.\cr\end{eqnarray*}

Under the change of variables $\lambda = [(n-1)\theta/\sqrt{2}]- (k-1)
\pi\sqrt{2}, \ k = 1,\ldots,n$, the
first and last integrals become $1/2$, and the integrals under the sum
become $1$.  Thus
\[L(f) =-\int_{S^1}(\id,f)^*\MQ_{\Delta_\epsilon}
=-{1\over 2}-\sum_{k=2}^{n-1}1 -{1\over 2}
=1-n,\]
as expected.

\bigskip

The formula for the Lefschetz number for functions on flat manifolds
is not much more complicated.  
Since $M$ is flat, for the Levi-Civita connection $\theta$ there
exists a local orthornomal frame $\{ x^i\}$ for which
the connection  and the curvature forms vanish.
The Mathai-Quillen form for the normal bundle to the diagonal
 is thus given by
\begin{equation}\label{flatmq}
\MQ_{\nu_\Delta}=\pi^{-n/2}e^{-x^2}dx^1\wedge ...\wedge
dx^n,\end{equation}
where $x$ is the fiber coordinate.

We need an explicit isomorphism 
$\alpha :\Delta_\epsilon\rightarrow\nu_\Delta$ between an 
$\epsilon$-neighborhood of the diagonal 
and the normal bundle to compute 
$\MQ_{\Delta_\epsilon}=\alpha^* \MQ_{\nu_\Delta}.$
Even though the exponential map is trivial near the diagonal, we will
use it initially to avoid confusion between normal vectors and points
of
$M\times M.$  

Fix $(x,y)\in \Delta_\epsilon$.  Since
the normal bundle consists of vectors of the form $(-v,v)$, there
exists $(\bar x,\bar x)\in\Delta$ such that $(x,y) = \exp_{(\bar
x,\bar x)}(-v,v) = (\exp_{\bar x}(-v),\exp_{\bar x}v).$  Thus $\bar x$
is the midpoint of the geodesic $\exp_{\bar x}(tv)$,
$t\in[-1,1]$  from $x$ to $y$.
This gives an isomorphism $\eta :U\to\Delta_\epsilon$
between a neighborhood $U$  of the zero section in $\nu_\Delta$ and
the 
$\epsilon$-neighborhood of the diagonal:
$$\eta(v
,-v)_{(\bar{x},\bar{x})}=(\exp_{\bar{x}}(v),\exp_{\bar{x}}(-v)).$$
Let $\rho:[0,\epsilon)\rightarrow [0,\infty)$  be a 
fixed diffeomorphism with $\rho(0)=0$, $\lim_{d\to\epsilon} 
\rho (d)=\infty$.
As before, we extend $\rho$
to take on values $\infty$ outside of $[0,\epsilon)$.

In the product
metric, $d((x,y),(\bar x,\bar x)) = d(x,y) /\sqrt{2}$, so
$\Delta_\epsilon = \{(x,y)\in M\times M: d(x,y) <\sqrt{2}\epsilon.\}.$
Thus $\beta: U\to\nu_\Delta$ given by
$$\beta(v,-v)_{(\bar{x},\bar{x})}=\Bigg\lbrace\matrix{
\Bigg(\rho\bigg({d(x,y)\over\sqrt{2}}\bigg){v\over |v|},
-\rho\bigg({d(x,y)\over\sqrt{2}}\bigg){v\over |v|}\Bigg),& v\neq 0,\cr
(\bar{x},\bar{x}),& v=0,\cr}$$
is a diffeomorphism and $\alpha =\beta\circ\eta^{-1}:\Delta_\epsilon
\to\nu_\Delta$ is our desired map:
$$\alpha(x,y)=\Bigg\lbrace\matrix{
\Bigg((\bar{x},\bar{x}),
\bigg( \rho\bigg({d(x,y)\over\sqrt{2}}\bigg){v\over |v|},
-\rho\bigg({d(x,y)\over\sqrt{2}}\bigg){v\over |v|}\bigg)\Bigg),& x\neq
y,\cr
((x,x),0),& x=y,\cr}$$
where
$x=\exp_{\bar{x}}(v),\ y=\exp_{\bar{x}}(-v),\ 
|v|=d(x,y)/\sqrt{2}.$

In the flat case, the map $\eta$ can be treated as
identity map, and there exists an isometry
between the $\epsilon$-tube $\Delta_\epsilon$ and $U$ such that the
map 
$\alpha$ reduces to $\beta:{\bf R}^n\to {\bf R}^n,$ with 
$$\beta(v)=\Bigg\lbrace\matrix{\rho(|v|){v\over |v|},& v\neq 0,\cr
0,&v=0,\cr}$$
and
$$\MQ_{\Delta_\epsilon} =\alpha^* \MQ_{\nu_\Delta}=\beta^*
\MQ_{\nu_\Delta}.$$

By (\ref{flatmq}), computing this last term reduces to  calculating
$\beta^* \dvol$,
which
we do in polar coordinates.
For 
$\gamma(t)=v +t{v\over |v|}$,
we have
\begin{eqnarray*}\beta_{* v}(\partial_r)
&=&{d\over dt}\Big|_{t=0}\beta(\gamma(t))
={d\over dt}\Big|_{t=0}\rho(|\gamma(t)|)
{\gamma(t)\over |\gamma(t)|}\cr
&=&\rho '(|v|)\Big({d\over dt}\Big|_{t=0}|\gamma(t)|\Big)
{\gamma(0)\over |\gamma(0)|}\cr
&\ &\ \ \ \ \ \ \ \ \ +\rho(|v|){\dot\gamma(0)\over |\gamma(0)|}+
\rho(|v|)\Big\lbrack\gamma(0)(-{1\over 2}
\langle\gamma(0),\gamma(0)\rangle^{-3/2}
2\langle\dot\gamma(0),\gamma(0)\rangle)\Big\rbrack\cr
&=&\rho '(|v|)\cdot 1\cdot {v\over |v|}+\rho(|v|){v\over |v|^2}
+\rho(|v|)\bigg({-v\over
|v|^3}\Big\langle{v\over|v|},v\Big\rangle\bigg)\cr
&=&\rho '(|v|)\partial_r+\rho(|v|)
\bigg\lbrack{v\over |v|^2}-{v\over |v|^2}\bigg\rbrack
=\rho '(|v|)\partial_r,\cr\end{eqnarray*}
and by a similar calculation using $\gamma(t)$ with
 $\dot\gamma(0)=(\partial_{\theta^i})_v$, we have
\[\beta_{* v}(\partial_{\theta^i})
=0+\rho(|v|){(\partial_{\theta^i})_v
\over|\partial_{\theta^i}|}+\rho(|v|){v\over |v|^3}\cdot 0
=\rho(|v|){\partial_{(\theta^i)_v}\over |v|}.\]
Thus
\begin{equation}\label{eqnone}\beta^* dr_v=\rho
'(|v|)dr_{\beta(v)}.\end{equation} 
Now
$(\partial_{\theta^i})_v=(\partial\theta^i\Big|_{v\over|v|})|v|$, as
$\langle \partial_{\theta^i}, \partial_{\theta^i}\rangle =r,$
and
$(\partial_{\theta^i})_{\beta(v)}=
(\partial\theta^i\Big|_{v\over|v|})\rho(|v|)$,
 as $|\beta(v)|=\rho(|v|).$
Hence
$(\partial_{\theta^i})_{\beta(v)}={\rho(|v|)\over |v|}
(\partial_{\theta^i})_v,$
so $(\beta_*)_v(\partial_{\theta^i})={\rho(|v|)\over |v|}
(\partial_{\theta^i})_v$ implies
$\beta_{*v}(\partial_{\theta^i})=(\partial_{\theta^i})_{\beta (v)}.$
Thus 
\begin{equation}\label{eqntwo}
\beta^* d\theta^i_v=d\theta^i_{\beta (v)}.\end{equation}
By (\ref{eqnone}), (\ref{eqntwo}), we get
\begin{eqnarray*}\MQ_{\Delta_\epsilon}
&=&\beta^* \MQ_{\nu_\Delta}\cr
&=&\beta^*(\pi^{-n/2}e^{-x^2}r^{n-1}dr
\wedge d\theta^1\wedge ...\wedge d\theta^{n-1})\cr
&=&\pi^{-n/2}e^{-\rho^2(|v|)}\rho(|v|)^{n-1}\rho'(|v|)
\dvol_{\alpha(x,y)}\cr
&=&\pi^{-n/2}e^{-\rho^2({d(x,y)\over\sqrt{2}})}\rho'\bigg({d(x,y)\over\sqrt{2}}\bigg)
\dvol_{\alpha(x,y)}.\cr\end{eqnarray*}

The last step is to calculate
\begin{eqnarray*}(\id,f)^* \MQ_{\Delta\epsilon}
&=&(\id,f)^* \Bigg\lbrack\pi^{-n/2}e^{-\rho^2({d(x,y)\over\sqrt{2}})}
\rho'\bigg({d(x,y)\over\sqrt{2}}\bigg)\dvol_{\alpha(x,y)}\Bigg\rbrack\cr
&=&\pi^{-n/2}e^{-\rho^2({d(x,f(x))\over\sqrt{2}})}
\rho'\bigg({d(x,f(x))\over\sqrt{2}}\bigg)(\id,f)^*(\dvol_{\alpha(x,y)}).
\cr\end{eqnarray*}
Here $\dvol_{\alpha (x,y)}$ is the volume element on the normal
bundle,
considered as a form near the diagonal. On a general manifold,
 calculating $(\id,f)^* \dvol_{\alpha(x,y)}$ will require introducing
coordinates on the tubular neighborhood
via the exponential map. Since $M$ is flat, the calculations
reduce to 
the case $M={\bf R}^n$.

For this let $(x^i)$ be flat coordinates near $x$, and let $(y^i)$ be
flat coordinates at $y$ given by parallel translating the
$\partial_{x^i}$ along the geodesic through $\bar x.$ (Here we are
assuming that $(x,f(x))$ is in the tubular neighborhood,
as $(\id,f)^*\MQ_{\Delta_\epsilon}$ vanishes otherwise.)
Then
$$\dvol_{\alpha(x,y)}=\bigwedge_{i=1}^n\bigg({-dx^i+dy^i\over\sqrt{2}}\bigg),$$
since the normal fiber $\nu_{(\bar{x},\bar{x})}$ at
 ${(\bar{x},\bar{x})}$ consists of vectors of the form  $(-v,v)$.
In the $(x^i),\ (y^i)$ coordinates, we may write 
 $f=(f^1,... ,f^n)$. Then
\begin{eqnarray*}(\id,f)^* \dvol_{\alpha(x,y)}
&=&(\id,f)^*\bigwedge_{i=1}^n\bigg({-dx^i+dy^i\over\sqrt{2}}\bigg)
=2^{-n/2}\bigwedge_{i=1}^n(-dx^i+df^i)\cr
&=&2^{-n/2}\bigwedge_{i=1}^n\bigg(-dx^i+{\partial f^i\over\partial
x^j}
dx^j\bigg)\cr
&=& 2^{-n/2}\bigwedge_{i=1}^n\bigg\lbrack
\bigg({\partial f^i\over\partial x^i}-1\bigg) dx^i+\sum_{i\neq j}
{\partial f^i\over\partial x^j} dx^j\bigg\rbrack.\cr\end{eqnarray*}

 Let $A$ be the matrix:
$$a_{ij}=\Bigg\lbrace\matrix{
(\partial f_i/\partial x^i)-1,& i=j,\cr
\partial f_i/\partial x^j,&i\neq j,\cr}$$
i.e.~$A=df\circ\Vert-\id$, where $\Vert$ denotes parallel translation
from $f(x)$ to $x$ along their geodesic. 
Then 
\begin{eqnarray*}(\id,f)^* \dvol_{\alpha(x,y)}
&=&2^{-n/2}\bigwedge_{i=1}^n\sum_{j=1}^n a_{ij}dx^j
=2^{-n/2}\det(A)\dvol_M\cr
&=&2^{-n/2}\det (df\circ\Vert-\id)\dvol_M,\cr\end{eqnarray*}
and so
\begin{eqnarray}\label{previous}(\id,f)^* \MQ_{\Delta_\epsilon}
&=&(\id,f)^*\alpha^* \MQ_{\nu_\Delta}\\
&=&(2\pi)^{-n/2}e^{-\rho^2({d(x,f(x))\over\sqrt{2}})}
  \rho'\Bigg({d(x,f(x))\over\sqrt{2}}\Bigg)\det\big(df\circ 
\Vert - \id\big)\dvol_M.\nonumber\end{eqnarray}

Since $(-1)^n\det(df\circ\Vert - \id) = \det(\id- df\circ\Vert)$,
Theorem \ref{theoremone} and (\ref{previous}) yield:
\begin{theorem}\label{theoremtwo}
Let $f:M\to M$ be a smooth
map of a closed, oriented, flat $n$-manifold $M$. Let $\rho
:
[0,\epsilon)\to [0,\infty)$ be an orientation preserving 
diffeomorphism and set $\rho (t)=\infty$ for $t\geq\epsilon$.
Then the Lefschetz number of $f$ is given by
\[L(f)
={1\over (2\pi)^{n/2}}\int_M
e^{-\rho^2\Big({d(x,f(x))\over\sqrt{2}}\Big)}\rho'
\bigg({d(x,f(x))\over\sqrt{2}}\bigg)
\det\big(\id -df\circ \Vert \big)\dvol_M.\]
\end{theorem}

Note that when $f$ is the identity map,
$\det(\id -df\circ \Vert)$ vanishes and we get
$\chi (M)=L(\id)=0,$
which reflects the fact that the Euler characteristic
of a flat manifold is zero.  Of course, the $\sqrt{2}$ factor in the
integrand can be incorporated into the diffeomorphism $\rho.$

We give some easy applications of our techniques. Let
$A_f$ be the portion of the graph of $f$  inside the tubular
neighborhood:
$$A_f=\{x\in M|(x,f(x))\in\Delta_\epsilon^M\}.$$
\begin{corollary}\label{corollaryone} (i) Let $f:M\to M$ be a smooth
map of a closed, oriented, flat manifold $M$ such that

a. For some $\epsilon >0$, $A_f$ is connected, and

b. $(\id-df\circ \Vert)$ is invertible on $A_f$.

\noindent Then $f$ has a fixed point.

(ii) If $g$ is $C^0$ close to the 
identity map, then there exists $x\in M$ and $v\in T_xM$, $v\neq 0$,
such that
$(dg\circ \Vert)(v)=v.$ \end{corollary}

\noindent {\sc Proof:}  (i) Since $\id -df\circ \Vert $ is invertible
on $A_f$,
$\sgn\ \det(\id -df\circ\Vert)$ is constant on the open set $A_f$.
The
integrand in Theorem \ref{theoremtwo} thus has constant sign on its
support, which is contained in $A_f.$  Thus $L(f) \neq 0$, and so
$f$ has a fixed point.

(ii)  Since $g$ is close to the identity,
$A_g=M.$
If $(dg\circ \Vert)(v)\neq v$ for all
$0\neq v\in TM$, then $\id -dg\circ \Vert $ is invertible.
This implies
$L(g)\neq 0,$ which contradicts
$L(g)=L(\id)=\chi(M)=0.$  \hfill$\Box$
\bigskip

Set $\|df\|=\sup_{x\in M}\|df_x\|.$  
\begin{proposition}\label{propositionone}
Let $f:M\to M$ be a smooth
map of a closed, oriented, flat $n$-manifold $M$. Then
$$|L(f)|\leq {C\over (2\pi)^{n/2}}\vol(M)(\| df\| +1)^n$$
for some constant $C>0$. \end{proposition}

\noindent{\sc Proof:} We may choose $\rho$ so that
$\lim_{z\to\epsilon}e^{-\rho^2(z)}\rho '(z)=0.$
Hence there exists $C>0$ such that
$0\leq e^{-\rho^2(z)}\rho '(z)\leq C.$
Note also that for  $v\in T_{f(x)}M$, 
$$|(\id -df\circ \Vert)(v)|\leq (\|df|\
+1)|v|,$$
since parallel translation
in an isometry. Thus
$|\det(\id - df\circ \Vert)|\leq (\| df\| +1)^n.$
By Theorem \ref{theoremtwo}, we have
\begin{eqnarray*}|L(f)|   &\leq& {1\over (2\pi)^{n/2}}\int_M\Big| 
e^{-\rho^2({d(x,f(x))\over\sqrt{2}})}
\rho '\Big({d(x,f(x))\over\sqrt{2}}\Big)\Big|\cdot\Big| 
\det(\id - df\circ \Vert)\Big| \dvol \cr
&\leq&{C\over (2\pi)^{n/2}}\vol(M) (\| df\| +1)^n.\cr
\end{eqnarray*}
\hfill$\Box$\\
\bigskip
The appendix  contains a similar result for arbitrary manifolds via Hodge
theory.

\subsection{\large\rm\bf Local expressions for arbitrary
metrics}\label{local}

In this subsection we calculate the local expression for the 
integrand in Theorem \ref{theoremone} for an arbitrary Riemanninan
metric.  This is more involved than in the flat case because the
exponential map is nontrivial.  The Jacobi fields which
measure the deviation of the exponential map from the identity enter
the computations.

A tubular  neighborhood $\Delta_\epsilon$ of the diagonal $\Delta$ in
$M\times M$ is
diffeomorphic to a neighborhood of the zero section in $\nu_M$, which
in turn is diffeomorphic to a neighborhood of zero in $TM$. The 
Levi-Civita connection on $M$ determines the space $H_M$ of
horizontal vectors on $TTM$,
while the space $V_M$ of vertical vectors is independent of the
connection. The
Mathai-Quillan form $\MQ_{TM}$ is written in terms of horizontal and
vertical vectors, so we have to identify the corresponding
horizontal and vertical vectors in the tube in order to compute
$\MQ_{\Delta_\epsilon}.$

Let $\alpha$ be the isomorphism from the neighborhood in $\nu_M$ to
the tube: for
$\nu_M
=\{ (v,-v):v\in TM\}$, we have  $\alpha
(v,-v)=(\exp_{\bar x} (v),\exp_{\bar x}(-v))$ at $(\bar x,\bar
x)\in\Delta$. 
As before, we take the radius of the tube small enough so that there
exists a unique minimal geodesic between $x$ and $y$ whenever $(x,y)$
is in the tube.
If $(x,f(x))\in \Delta_\epsilon$, we know that $\bar
x$ is the midpoint  of the unique minimal geodesic $\gamma$ from $x$
to $f(x)$ and
$|v|=d(x,f(x))/{2}$.

Pick an orthonormal frame $\{ Y_i\}$ at $\bar x$. Let $\beta
:TM\rightarrow\nu_M$ be the bundle isomorphism $\beta
(v_x)=(v_x,-v_x)$. The horizontal space $H$ in the tube is defined to
be $d(\alpha\beta )(H_M)$,
and the vertical space $V$ in the tube is $d(\alpha\beta
)(V_M)$. Define 
vectors $X_i,\widetilde X_i$ at $x$ by
\begin{eqnarray}\label{oneone}
X_i &=& d(\exp_{\bar x} )_v(Y_i), \nonumber \\
\widetilde{X}_i &=& d(\exp_\cdot \Vert v)_{\bar x}(Y_i), 
\end{eqnarray}
where in the first line $Y_i$ is trivially translated to a vector in
$T_vT_{\bar x}M$, and in the second line $\Vert v$ denotes the
parallel translation of $v$ along a curve in $M$ with tangent vector
$Y_i$. Similarly define vectors $Z_i,\widetilde{Z}_i$ at $f(x)$ by
replacing $v$ in (\ref{oneone})
with $-v$. If we parametrize $\gamma$ from $\bar x$ to
$x$ as $\gamma (t)$, then $\widetilde{X}_i$ is the endpoint of a
Jacobi field $J$ with $J(0)=Y_i$ -- i.e. $J$ is the variation vector
field of the family of geodesics $\gamma_s(t) = 
\exp_{\eta (s)} (t\Vert v)$, where
$\dot\eta (0)=Y_i$ and $t\in [0,1]$. Similarly, $X_i$ is the
endpoint of a Jacobi field $J$, the variation vector field of the
family of geodesics $\gamma_s(t)=\exp_{\bar x}(t(v +sY_i))$, which has
$(\nabla
J)(0)=Y_i$ (cf. \cite[Cor.~3.46]{Ga}). Similar remarks apply to
$Z_i,\widetilde{Z}_i$.
\begin{lemma}
The vertical space $V$ at $(x,f(x))$ is spanned by $\{ (-X_i,Z_i)\}$
and
the horizontal space $H$ is spanned by $\{
(\widetilde{X}_i,\widetilde{Z}_i)\}$.
\end{lemma}

\noindent{\sc Proof:} Set $\delta =\alpha\beta$. A vertical vector at
$v\in T_{\bar x}M$ is a tangent vector $Y$ to a curve $\eta
(t)\subset T_xM$ with $\eta (0)=v,\dot\eta (0)=Y$. Then
\begin{eqnarray*}
d\delta_v(Y) &=&\frac{d}{dt}\biggl|_{t=0} (\exp_{\bar x}\eta (t),\
\exp_{\bar x}(-\eta (t))) \\
&=& \left( \frac{d}{dt}\biggl|_{t=0} \exp_{\bar
x}\eta(t),\frac{d}{dt}\biggl|_{t=0} \exp_{\bar x}(-\eta (t))\right)
\\
&=& (d(\exp_{\bar x})_v Y,\ d(\exp_{\bar x})_v (-Y)).
\end{eqnarray*}
Thus the vertical space at $(x,f(x))=(\exp_{\bar x}v,\exp_{\bar
x}(-v))$ is spanned by $\{ (d(\exp_{\bar x})_v(Y_i), d(\exp_{\bar
x})_{-v}(Y_i))\}$.

Let $\Vert v=\Vert_y v$ denote the parallel translation of $v$ along
radial geodesics centered at $\bar x$. Then $\Vert v$ is parallel at
$\bar x$, and the horizontal vectors at $v$ are spanned by
$$\frac{d}{dt}\biggl|_{t=0} \Vert_{ \exp_{\bar x}(tY_i)}v.
$$
Thus the horizontal vectors at $(x,f(x))$ are spanned by
\begin{eqnarray*}
\frac{d}{dt}\biggl|_{t=0} \delta (\Vert_{ \exp_{\bar
x}(tY_i)}v) &=& \frac{d}{dt}\biggl|_{t=0} (\exp_{\exp_{\bar
x}(tY_i)}\Vert_{
\exp_{\bar x}(tY_i)}v,\exp_{\exp_{\bar x}(tY_i)}\Vert_{\exp_{\bar
x}(tY_i)}(-v)) \\
&=& (\widetilde{X}_i,\widetilde{Z}_i).\end{eqnarray*}
${}$\hfill$\Box$\\

\noindent {\bf Remarks:} 1) The lemma shows that at $(x,f(x))$,
\begin{eqnarray*}
V &=&(-d(\exp_{\bar x})_v, d(\exp_{\bar x})_{-v})V_M, \\
H&=& (d(\exp \Vert v)_{\bar x}, d(\exp -\Vert v)_{\bar x})H_M.
\end{eqnarray*}

2) It is easy to check that
$X_i,\widetilde{X}_i,Z_i,\widetilde{Z}_i$ are just parallel
translations of $Y_i$ if $M$ is flat.

3) Vertical vectors at $(x,f(x))$ are those pairs of vectors in
$T_xM\times T_{f(x)}M$
which are endpoints of a Jacobi field along $\gamma$ which
vanishes at $\bar x$. Horizontal vectors are pairs
of vectors which are endpoints of a Jacobi field along $\gamma$ which
is parallel at $\bar x$.
\bigskip

\begin{lemma}

Let $(X,Z)\in T_{(x,f(x))}(M\times M)$. Take the unique Jacobi
field $Y$ along $\gamma$ with $Y(x)=X$, $Y(f(x))=Z$. Let $X_1,Z_1$ be
the values at $x,f(x)$ of the Jacobi field   $Y_1$ along $\gamma$
given by $Y_1(\bar x)=0$, $\frac{DY_1}{dt}(\bar x)=\frac{DY}{dt}(\bar
x)$. Let $\widetilde X_1,\widetilde{Z}_1$ be the values at $x,f(x)$ of
the Jacobi
field $Y_2$ along $\gamma$ given by $Y_2(\bar x)=Y(\bar x)$,
$\frac{DY_2}{dt}(\bar x)=0$. Then
$(X,Z)=(X_1,Z_1)+(\widetilde{X}_1,\widetilde{Z}_1)$ is the
decomposition of $(X,Z)$ into vertical and horizontal vectors.
\end{lemma}

\noindent {\sc Proof:} By Remark 3,
$(X_1,Z_1),(\widetilde{X}_1,\widetilde{Z}_1)$ are vertical and
horizontal vectors respectively. Since the Jacobi equation is linear,
the endpoints of the Jacobi field $Y_1+Y_2$ are $X_1+\widetilde{X}_1,
Z_1+\widetilde{Z}_1$. Since $Y_1+Y_2$ has the same position and
velocity vectors as $Y$ as $\bar x$, we must have
$X=X_1+\widetilde{X}_1, Z=Z_1+\widetilde{Z}_1$. \hfill$\Box$\\

Let $\rho :[0,\epsilon )\rightarrow [0,\infty )$ be a diffeomorphism
fixing zero,
and which extends smoothly  to an even 
function on $(-\epsilon,\epsilon )$. Let $\MQ_\nu$
be the Mathai-Quillen form of the normal bundle $\nu = \nu_\Delta$ and let
$\MQ_{\Delta_\epsilon}
=(\exp^{-1})^\ast \rho^\ast \MQ_\nu$ be the corresponding
Mathai-Quillen form on $M\times M$. Here we abbreviate
$(\exp^{-1}, \exp^{-1})$ to just
$\exp^{-1}$. 

In (\ref{basicdef}), the vertical
coordinates are denoted by $x^i$ and the horizontal coordinates are
hidden in $\Omega_I.$  For the calculations on $TM$, we need to make
the horizontal coordinates explicit and take care not to confuse them
with the vertical coordinates.
So let $\{ x^i\}$ be a synchronous orthonormal frame
centered at $\bar x$. In each fiber of $\nu$, we take the orthonormal
polar coordinate frame  $\{(-x^i,x^i)\}$, with
$\{ x^i\} =\{ x^1
=\partial_r, r^{-1}\partial_{\theta^i}\}$,
away from the origin. These frames do not
agree at the origin in each fiber, but the formulas below will be
smooth at the origin. We know $(\rho^\ast dx^i)_{v=0}=0$ since
$\rho^\prime (0)=0$, and for $v\not= 0$,
\begin{eqnarray*}
\rho^\ast dx_v^1 &=& \rho^\ast dr_v =\rho^\prime (|v|)dr_{\rho (v)},\\
\rho^\ast d\theta_v^i &=& d\theta^i_{\rho (v)}.\end{eqnarray*}
Thus 
$$[(\exp^{-1})^\ast \rho^\ast dx^i](\alpha
)=
\left\{ \begin{array}{ll}
\rho^\prime (|(\exp^{-1})_\ast\alpha
|)(\exp^{-1})^\ast dx^i(\alpha), &\mbox{if $i =1$,}
\\
\hfil \\
(\exp^{-1})^\ast dx^i(\alpha), &\mbox{if $i \not= 1$}.
\end{array}\right.$$
If $\{
y^i\}$ is another synchronous frame centered at $\bar x$ (possibly
equal to $\{ x^i\}$), then the horizontal lifts of $y^i$ into $T\nu$
are orthonormal in the metric on $T\nu$ induced by the metric on
$M$. Since $\rho_\ast (y^i)=y^i$, we have
$$[(\exp^{-1})^\ast \rho^\ast dy^i](\alpha
)=(\widetilde{X}_i,\widetilde{Z}_i)^{\#} (\alpha ),
$$
where $(\widetilde{X}_i,\widetilde{Z}_i)^{\#}$ is the cotangent vector
dual to $(\widetilde{X}_i,\widetilde{Z}_i)$.
If $({\rm Id},f)^\ast (\MQ_{\Delta_\epsilon} )_{(x,x)}=D\ \mbox{\rm
dvol}_{(x,x)}$, then
$$D=({\rm Id},f)^\ast (\MQ_{\Delta_\epsilon} )_{(x,x)}\left(
\frac{(y^1,
y^1)}{\sqrt 2},\ldots,\frac{(y^n,y^n)}{\sqrt 2}\right).
$$
Thus
\begin{eqnarray*}
D
&=& 
\frac{e^{-\rho^2(\frac{d(x,f(x))}{\sqrt{
2}})}}{(2\pi)^{n/2}}\sum_{I,I^\prime}
\epsilon (I,I^\prime )[(\exp_{\bar x}^{-1})^\ast \rho^\ast
]{\rm Pf}(
\Omega_I)\\
&&\qquad \wedge dx^{I^\prime}((y^1,f_\ast y^1 ),\ldots,
(y^n,f_\ast y^n)).\end{eqnarray*}
Write 
\begin{equation}\label{alpha}(y^i,f_\ast
y^i)=P_V^i+P_H^i\end{equation}
 for the decomposition of
$(y^i,f_\ast y^i)$ into vertical and horizontal vectors as in the
lemma.  Let $\Sigma_n$ be the permutation group on $\{1,\ldots,n\}.$
Then
\begin{eqnarray*}
D &=&
\ALpha \sum_{I,I^\prime}
\frac{\epsilon (I,I^\prime)}{|I|!\ |I^\prime|!}
\sum_{\mu\in \Sigma_n}(\exp_{\bar x}^{-1})^\ast {\rm
Pf}(
\Omega_I)(P_H^{\mu_1},\ldots,P_H^{\mu_{|I|}})\\
&&\qquad \cdot \widetilde{\rho}^\prime_{I^\prime} \left(
\frac{d(x,f(x))}{\sqrt
2}\right) d\widetilde{x}^{I^\prime} \left( P_V^{\mu_{|I|+1}},\ldots,
P_V^{\mu_n}\right),\end{eqnarray*}
where
$$\widetilde{\rho}^\prime_{I^\prime} \left( \frac{d(x,f(x))}{\sqrt
2}\right)
=
\left\{ \begin{array}{ll}
\rho^\prime (\frac{d(x,f(x))}{\sqrt 2}),
 &\mbox{if $i_1^\prime =1$,}
\\
\hfil \\
1, &\mbox{if $i_1^\prime \not= 1$}.
\end{array}\right.
$$
Here $d\widetilde{x}^{I^\prime}=(\exp_{\bar x}^{-1})^\ast dx^i$, and
we have used $\rho^\ast$ {\rm Pf}$(
\Omega_I)=$ {\rm Pf}$
(
\Omega_I)$,
since this
Pfaffian is a horizontal form. 

Define $n\times n$ matrices $A=A_x$,
$B=B_x$ by
\begin{equation}\label{beta}(\exp_{\bar x}^{-1})_\ast P_H^i
=A_j^iy^j,\ \ (\exp_{\bar
x}^{-1})_\ast P_V^i=B_j^i x^j,\end{equation}
where strictly speaking the last term is $(-B_j^i x^j,B_j^i x^j).$
For example, at a fixed
point $x=f(x)$, the decomposition of $(q,f_\ast q)$ into vertical and
horizontal components is given by
$$(q,f_\ast q)=\left( \frac{q-f_\ast q}{2},\frac{-q+f_\ast
q}{2}\right) +\left( \frac{q+f_\ast q}{2},\frac{q+f_\ast q}{2}\right),
$$
since
$(\exp^{-1})_\ast ={\rm Id}$ at a fixed point.  Since vertical
(resp. horizontal) vectors on the diagonal are of the form 
$(-v,v)$ (resp. $(v,v)$),  $A$ is the
matrix of
$\frac{1}{2} (df+{\rm Id})$ and $B$ is the matrix of
$\frac{1}{2}(df-{\rm Id})$. It follows easily from Remark 2
 that on a flat manifold,
$A=\frac{1}{2} (\Vert\circ df+{\rm Id}),B=\frac{1}{2}(\Vert \circ
df-{\rm Id})$
for arbitrary $x$.

Then
\begin{eqnarray}\label{onetwo}D &=&
\ALpha \sum_{I,I^\prime}c_{|I|}\rats
\sum_{\mu\in\Sigma_n}
({\rm sgn}\ \mu)A^{\mu(1)}_{j_1}\cdot\ldots\cdot
A^{\mu(|I|)}_{j_{|I|}}\nonumber\\
&&\qquad\cdot
{\rm Pf} (
\Omega_I)_{\bar x}(y^{j_1},\ldots,
y^{j_{|I|}})\nonumber\\
&&\qquad \cdot
\widetilde{\rho}^\prime_{I^\prime} \left( \frac{d(x,f(x))}{\sqrt
2}\right)
B_{k_1}^{\mu(|I|+1)}\cdot\ldots\cdot
B_{k_{|I^\prime|}}^{\mu(n)}dx^{I^\prime}
(x^{k_1},\ldots, x^{k_{|I^\prime |}}).
\end{eqnarray}
(We use summation convention for the $j$ and $k$ indices.)
The right hand side of (\ref{onetwo}) vanishes unless
$I^\prime = \{k_1,\ldots,k_{|I^\prime|}\}\equiv K$, and
$$\sum_{K, K=I^\prime}B_{k_1}^{\mu(|I|+1)}\cdot\ldots\cdot
B_{k_{|I^\prime|}}^{\mu(n)}dx^{I^\prime}
(x^{k_1},\ldots, x^{k_{|I^\prime |}}) =
\det(B^{\mu(|I|+q)}_{i^\prime_s}),$$
for $q,s = 1,\ldots, |I^\prime|.$  Here $I^\prime =
\{i^\prime_1,\ldots,i^\prime_{|I^\prime|}\}$, with $i^\prime_1 <\ldots
< i^\prime_{|I^\prime}.$  We denote this determinant by
$\det(B^\mu_{I^\prime}).$
For $I= \{i_1,\ldots,i_{|I|}\}$, with $i_1 <\ldots < i_{|I|}$, we have
\begin{eqnarray*}
{\rm
Pf}(
\Omega_I)&=& c_{|I|}\sum_{\sigma,\tau\in\sum_{|I|}}(\mbox{\rm
sgn}\ \sigma )(\mbox{\rm sgn}\ \tau)
R_{i_{\sigma(1)}i_{\sigma(2)}i_{\tau(1)}i_{\tau(2)}}\cdot\ldots\cdot
R_{i_{\sigma(|I|-1)}i_{\sigma(|I|)}i_{
\tau(|I|-1)}i_{\tau(|I|)}}\\
&&\qquad \cdot dy^{1}\land
\ldots\land dy^{{|I|}},
\end{eqnarray*}
with $c_{|I|} = (-1)^{|I|/2}[2^{|I|}(|I|/2)!]^{-1}$ \cite[(1.3)]{MQ}.
(The $(-1)^{|I|/2}$ reflects our sign convention on
curvature.)  As above, the right hand side of (\ref{onetwo}) vanishes
unless $I= \{j_1,\ldots,j_{|I|}\}\equiv J.$   Summing over such $J$
produces the term $\det(A^{\mu(t)}_{i_k})\equiv \det(A^\mu_I)$, 
for $t,k = 1,\ldots,|I|$.

Thus
\begin{eqnarray*}
D
&=& 
\ALpha \sum_{I,I^\prime}
 \sum_{\mu\in \Sigma_n} c_{|I|}\rats
({\rm sgn}\ \mu)
\det (A_{I}^{\mu})\det
(B_{I^\prime}^{\mu})
\widetilde{\rho}^\prime_{I^\prime}
\left( \frac{d(x,f(x))}{\sqrt
2}\right)\\
&&\qquad\cdot \sum_{\sigma,\tau\in\Sigma_{|I|}}(\mbox{\rm sgn}\ \sigma
)(\mbox{\rm sgn}\  \tau)
R_{i_{\sigma(1)}i_{\sigma(2)}i_{\tau(1)}i_{\tau(2)}}\cdot\ldots\cdot
R_{i_{\sigma(|I|-1)}i_{\sigma(|I|)}i_{ \tau(|I|-1)}i_{\tau(|I|)}}.
\end{eqnarray*}
Since $\int_M dy^1 \land\ldots\land
dy^n=2^{-n/2}\int_\Delta d(y^1,y^1
)\land\ldots\land d(y^n,y^n)$, we have
\begin{theorem}\label{localtheorem}
 Let $f:M\to M$ and fix $\epsilon >0$ such that $(x,y) \in
\Delta_M^\epsilon$ implies the existence of a unique minimal geodesic
between $x$ and $y$.
For $x\in M$,
define matrices $A,\ B$ by (\ref{alpha}), (\ref{beta}), provided
$(x,f(x))$ is in the $\epsilon$-neighborhood of the
diagonal; otherwise set $A=B=0.$
Then 
\begin{eqnarray*}
L(f) &=& \frac{(-1)^{{\rm dim}\ M}}{(2\pi)^{n/2}}\int_M
e^{-{\rho^2(\frac{d(x,f(x))}{\sqrt 2})}}
\sum_{I,I^\prime} c_{|I|}
\rats\sum_{\mu\in\Sigma_n} ({\rm sgn}\ \mu)
\det (A_{I}^{\mu})\det
(B_{I^\prime}^{\mu})\\
&&\qquad \cdot\widetilde{\rho}^\prime_{I^\prime}
\left( \frac{d(x,f(x))}{\sqrt
2}\right)
\cdot\sum_{\sigma,\tau\in\Sigma_{|I|}}
R_{i_{\sigma(1)}i_{\sigma(2)}i_{\tau(1)}i_{\tau(2)}}\cdot\ldots\cdot
R_{i_{\sigma(|I|-1)}i_{\sigma(|I|)}i_{ \tau(|I|-1)}i_{\tau(|I|)}}\\
 &&\qquad \cdot \dvol.
\end{eqnarray*}
\end{theorem}

\noindent  The hypothesis on $\epsilon$ is satisfied
once $\epsilon$ is less than half the injectivity
radius of $M$.
If the
injectivity
radius $i$ is known, we can set $\rho (x)={\rm sec}(
\pi x/2i) -1$, for example.

As a check, we consider the case where $f={\rm Id}$. Then $A = {\rm
Id}$ and
$B=0$, so if $n = {\rm dim}\ M$ is odd, the integrand vanishes
and we get $\chi(M) = L({\rm Id}) =0.$  If $n$ is even, the only
contribution to the integrand occurs when $I = \{1,\ldots,n\}, \
I^\prime =\emptyset$.  For this $I$, we
have $({\rm sgn}\ \mu) \det(A^\mu_I) = 1.$
The theorem  becomes
\begin{eqnarray*}
\chi(M)=L({\rm Id}) &=& \frac{(-1)^{n/2}}{(8\pi)^{n/2} (n/2)! }\int_M
\sum_{\tau,\sigma\in\Sigma_n} ({\rm sgn}\ \sigma)({\rm
sgn}\ \tau)
e^0
R_{i_{\sigma(1)}i_{\sigma(2)}i_{\tau(1)}i_{\tau(2)}}\cdot\ldots\\
&&\qquad \cdot
R_{i_{\sigma(|I|-1)}i_{\sigma(|I|)}i_{
\tau(|I|-1)}i_{\tau(|I|)}}\dvol,
\end{eqnarray*}
the Chern-Gauss-Bonnet theorem.  

The other extremal case occurs when $M$ is flat. Because $R_{ijkl}=0$,
the only contribution to the integrand occurs when $I=\emptyset$. Then
$c_0 = 1$ and $({\rm sgn}\ \mu ) \det (B^{\mu}_{I^\prime}) = \det
(\Vert\circ df - {\rm Id}) $, so
the integrand becomes
$$\frac{1}{(2\pi)^{n/2}}
e^{-{\rho^2(\frac{d(x,f(x))}{\sqrt 2})}}
\rho^\prime
\left(
\frac{d(x,f(x))}{\sqrt 2}\right) \det (\id -\Vert\circ
df),
$$
which agrees with the flat case formula, since parallel translation is
an isometry.

Note that at a fixed point, the integrand becomes
$$\sum_{\scriptstyle{I,I^\prime\atop 1\not\in I^\prime}}
c_{|I|} \rats \sum_{\mu\in\Sigma_n} ({\rm sgn}\ \mu)\det\biggl(
\frac{1}{2}(df-{\rm Id})_{I^\prime}^\mu \biggr)
\det\biggl(
\frac{1}{2}(df+{\rm Id})_{I}^\mu \biggr)\ \ \ \ \ $$
$$\cdot \rho^\prime
\left(
\frac{d(x,f(x))}{\sqrt 2}\right)
\sum_{\tau,\sigma\in\Sigma_{|I|}} ({\rm sgn}\ \sigma)({\rm
sgn}\ \tau)
R_{i_{\sigma(1)}i_{\sigma(2)}i_{\tau(1)}i_{\tau(2)}}\cdot\ldots\cdot
R_{i_{\sigma(|I|-1)}i_{\sigma(|I|)}i_{ \tau(|I|-1)}i_{\tau(|I|)}}.$$

\bigskip
\noindent {\bf Example:}  We determine the integrand in the Lefschetz
formula for
$M$ an oriented surface of constant sectional/Gaussian
curvature $-1$.   At the end we indicate the changes for constant
curvature manifolds in general.

We first determine the horizontal and vertical
components of a vector $(q,f_\ast q)\in T_{(x,f(x))}(M\times
M)$. Assume there exists a unique minimal geodesic $\gamma$ joining
$x$ to $f(x)$ with midpoint $\bar x.$ Let
$|\dot\gamma|=1$, and let $\alpha$ be the unit  normal to $\gamma$
determined by the orientation. Set $d=d(x,f(x))$.

Let $J$ be a Jacobi field along $\gamma$. Plugging
$J(t)=a(t)\dot\gamma +b(t)\alpha$ into the Jacobi equation
$D^2J/dt^2+R(\dot\gamma, J)\dot\gamma =0$ and using $\langle
R(\dot\gamma, J)\dot\gamma, \dot\gamma\rangle =0$, $\langle
R(\dot\gamma,\alpha )\dot\gamma,\alpha\rangle =-1$ yields $\ddot a=0,\
\ddot b-b=0$. Thus $J(t)=(c_0+c_1t)\dot\gamma +(d_1\hbox{\rm
sinh}(t)+d_2\cosh (t))\alpha$. Imposing the boundary conditions
$J(0)=q$, $J(d)=f_\ast q$ gives
$$J(t)=\left( q_1+\left( \frac{w_1-q_1}{d}\right)t\right)\dot \gamma
+\left[ \left( \frac{w_2-q_2\cosh (d)}{\sinh (d)}\right)\sinh
(t)+ q_2\cosh (t)\right]\alpha,
$$
where $q=q_1\dot\gamma +q_2\alpha$, $f_\ast q=w_1\dot\gamma
+w_2\alpha$. In particular
\begin{eqnarray*}
J\left( \frac{d}{2}\right) &=&\left(
\frac{q_1+w_1}{2}\right)\dot\gamma + \left[\left( \frac{w_2-q_2\cosh
(d)}{\sinh (d)}\right)\sinh \left( \frac{d}{2}\right) +q_2\cosh \left(
\frac{d}{2}\right)\right]\alpha \\
&=& \left( \frac{q_1+w_1}{2}\right)\dot\gamma +\left(
\frac{w_2+q_2}{2\cosh (\frac{d}{2})}\right)\alpha, \\
\frac{DJ}{dt}\left( \frac{d}{2}\right)
&=& \left(
\frac{w_1-q_1}{d}\right)\dot\gamma +\left( \frac{w_2-q_2}{2\sinh
(\frac{d}{2})}\right) \alpha,
\end{eqnarray*}
The Jacobi fields $J_1,J_2$ determined by $J_1(d/2)=0,\ (DJ_1/dt)
(d/2)=(DJ/dt) (d/2)$ and $J_2(d/2)=J(d/2),\ (DJ_2/dt)(d/2)=0$ are
given by
\begin{eqnarray*}
J_1(s) &=&\left( \frac{w_1-q_1}{d}\right) s\dot\gamma +\left(
\frac{w_2-q_2}{2\sinh (d/2)}\right)\sinh (s)\alpha, \\
J_2(s) &=& \left( \frac{q_1+w_1}{2}\right)\dot\gamma +\left(
\frac{w_2+q_2}{2\cosh (d/2)}\right) \cosh (s)\alpha,
\end{eqnarray*}
where $s=0$ corresponds to $\bar x$. Evaluating $J_1,J_2$ at $s=\pm
d/2$ gives the decomposition of $(q,f_\ast q)$ into vertical and
horizontal components:
\begin{eqnarray*}
(q,f_\ast q)_{\small\rm vert} &=&\left( \left(
\frac{-w_1+q_1}{2}\right)\dot\gamma -\left(
\frac{w_2-q_2}{2}\right)\alpha, \left(
\frac{w_1-q_1}{2}\right)\dot\gamma +\left(
\frac{w_2-q_2}{2}\right)\alpha\right), \\
(q, f_\ast q)_{\small\rm hor} &=& \left( \left(
\frac{q_1+w_1}{2}\right)\dot\gamma +\left(
\frac{w_2+q_2}{2}\right)\alpha, \left(
\frac{q_1+w_1}{2}\right)\dot\gamma +\left( \frac{w_2+q_2}{2}\right)
\alpha \right).
\end{eqnarray*}

Let $(x,f(x)) = (\exp_{\bar x}v, \exp_{\bar x} (-v))$, so $v =
-(d/2)\dot\gamma$ at $\bar x.$
We now determine $(\exp_{\bar x}^{-1})_\ast^{(2)}$:
$T_{(x,f(x))}(M\times M)\rightarrow T_{(v,-v)}
T_{(\bar x,\bar x)}(M\times M)$,
where $(\exp_{\bar x}^{-1})_\ast^{(2)}$ is shorthand for $(-(\exp_{\bar
x}^{-1})_\ast, (\exp_{\bar x}^{-1})_\ast )$. 
For a
vertical vector $\beta =\beta_1\dot\gamma +\beta_2\alpha\in T_vT_{\bar
x}M$, with $\dot\gamma,\alpha$ trivially parallel translated to $v$,
$$(\exp_{\bar x})_{\ast, v}(\beta )=\frac{d}{ds}\biggl|_{s=0}
\exp_{\bar x}(v + s\beta ),
$$
which is the value at $x$ of the Jacobi field $J$ along $\gamma$ with
$J(\bar x)=0$, $(DJ/dt)(\bar x)=2\beta /d$, since $|v|=d/2$. Solving
for $J$ as above, we get
$$(\exp_{\bar x})_{\ast, v}(\beta )=\beta_1\dot\gamma
+\frac{2}{d}\beta_2\sinh \left( \frac{d}{2}\right)\alpha.
$$
Thus
$$
(\exp_{\bar x}^{-1})_\ast^{(2)}(q,f_\ast q)_{\small\rm vert}=\ \ \ \
\ \ \ \ \ \ \ \ \ \ \ \ \ \ \ \ \ \ \ \ \ \ \ \ \ \ \ \ \ \ \ \ \ \ \
\ \ \ \ \ \ \ \ \ \ \ \ \ \ \ \ \ \ \ \ \ \ \ \ \ \ \ \ \ \ \ \ \ \ \
$$
$$\ \ \ \ \ \ \ \ \ \ \ \ \ \ \ \ \left(
\left( \frac{-w_1+q_1}{2}\right)\dot\gamma -\frac{d(w_2-q_2)}{4\sinh
(\frac{d}{2})} \alpha, \left( \frac{-w_1+q_1}{2}\right)\dot\gamma
-\frac{d(w_2-q_2)}{r\sinh (\frac{d}{2})}\alpha\right).
$$
Similarly, for a horizontal vector $\delta =\delta_1\dot\gamma
+\delta_2\alpha$, where $\dot\gamma,\alpha$ now denote the horizontal
lifts of $\dot\gamma,\alpha$ to $T_vT_{\bar x}M$, we have
$$(\exp_{\bar x})_{\ast ,v}(\delta )=\delta_1\dot\gamma
+\delta_2\cosh\left( \frac{d}{2}\right)\alpha,
$$
so
$$(\exp_{\bar x}^{-1})_\ast^{(2)} (q,f_\ast q)_{\small\rm hor} =\ \ \
\
\ \ \ \ \ \ \ \ \ \ \ \ \ \ \ \ \ \ \ \ \ \ \ \ \ \ \ \ \ \ \ \ \ \ \
\ \ \ \ \ \ \ \ \ \ \ \ \ \ \ \ \ \ \ \ \ \ \ \ \ \ \ \ \ \ \ \ \ \ \
$$
$$\ \ \ \ \ \ \ \ \ \ \ \ \ \ \ \ 
\left(
\left( \frac{q_1+w_1}{2}\right)\dot\gamma +\left(
\frac{w_2+q_2}{2\cosh (\frac{d}{2})}\right)\alpha, \left(
\frac{q_1+w_1}{2}\right) \dot\gamma +\left( \frac{w_2+q_2}{2\cosh
(\frac{d}{2})}\right) \alpha \right).
$$

To determine the matrix $B$, we have to express $(\exp_{\bar
x}^{-1})_\ast^{(2)} (q,f_\ast,q)$, for $q=\dot\gamma,\alpha$, in polar
coordinates at $(v,-v)$ in $\nu_{(\bar x,\bar x)}$. The radial vector
at $(v,-v)$ is
$$r=\left( \frac{v}{|v|\sqrt 2},\frac{-v}{|v|\sqrt 2}\right) =\left(
\frac{v\sqrt 2}{d},\frac{v\sqrt 2}{d}\right),
$$
and the unit angular vector ``$r^{-1}\partial_\theta$'' is $(-\alpha
/\sqrt 2,\alpha /\sqrt 2)$. Note that $(\dot\gamma, -\dot\gamma
)=-\sqrt 2 r$. For  $f_\ast \dot\gamma =w_{11}\dot\gamma
+w_{12}\alpha, f_\ast \alpha =w_{21}\dot\gamma +w_{22}\alpha$, we get
$$(B_j^i) =\left( \matrix{ -\frac{\sqrt 2}{2}(-w_{11}+1) & \frac{\sqrt
2d w_{12}}{4\sinh (d/2)} \cr
\hfil \cr
\frac{\sqrt 2}{2}w_{21} & \frac{\sqrt 2 d(w_{22}-1)}{4\sinh (d/2)}\cr}
\right).
$$
Thus
\begin{eqnarray*}
\det B &=& \frac{d}{4\sinh
(\frac{d}{2})}((w_{11}-1)(w_{22}-1)-w_{12}w_{21}) \\
&=& \frac{d}{\sinh (\frac{d}{2})}\det \left( \frac{1}{2}(\Vert \circ
df-{\rm Id})\right).
\end{eqnarray*}
Similarly,
$$(A_j^i) =\left( \matrix{ -\frac{w_{11}+1}{\sqrt 2} & \frac{-\sqrt 2
w_{12}}{2\cosh (d/2)} \cr
\frac{-w_{21}}{\sqrt 2} & \frac{ -\sqrt 2 (w_{22}+1)}{2\cosh
(\frac{d}{2})}\cr}
\right),
$$
and
$$\det A=\frac{2}{\cosh (\frac{d}{2})}\det \left( \frac{1}{2} (\Vert
\circ df+{\rm Id})\right).
$$
We now plug this information into the Lefschetz formula. Note that
$I=\{ 1,2\}$ or $I=\emptyset$ and that $R_{1212}=1$ in our
convention.
We obtain
\begin{eqnarray*}
L(f) &=& \frac{1}{2\pi}\int_M e^{-\rho^2(\frac{d}{\sqrt{2}})}\left[
\left(
\frac{2}{\cosh (\frac{d}{2})}\right) \frac{(-1)\cdot 4}{4\cdot 2!}\det
\left( \frac{1}{2} (\Vert \circ df+{\rm Id})\right) \right. \\
&&\qquad +\left. \rho^\prime
 \left( \frac{d}{\sqrt 2}\right)\left( \frac{d}{\sinh
(\frac{d}{2})}\right) \cdot \frac{1}{2!}\det\left( \frac{1}{2}
(\Vert\circ df -{\rm Id})\right)\right]dA.
\end{eqnarray*}
In the first line, there are factors of $c_{\{ 1,2\}}=-1/4, \ |I|!=2$,
and $\sum_{\sigma,\tau\in\Sigma_2} R_{\sigma (1)\sigma (2)\tau (1)\tau
(2)}=4$. In the second line, $c_\emptyset =1$ and $|I^\prime
|!=2$. Thus we obtain
\begin{proposition} \label{cons}
Let $M$ be an oriented surface of constant
curvature $-1$. Then
\begin{eqnarray*}
L(f) &=& \frac{1}{2\pi} \int_M e^{-\rho^2 (\frac{d(x,f(x))}{\sqrt
2})}\left[ \frac{-\det (\frac{1}{2}(\Vert\circ df_x+{\rm Id}))}{\cosh
(\frac{d(x,f(x))}{2})}\right. \\
&&\qquad +\left. \rho^\prime \left( \frac{d(x,f(x))}{\sqrt
2}\right)\left(
\frac{d(x,f(x))}{2\sinh (\frac{d(x,f(x))}{2})}\right)\det\left(
\frac{1}{2} (\Vert\circ df_x-{\rm Id})\right)\right]dA.
\end{eqnarray*}
\end{proposition}

It is straightforward to extend this result to higher dimensional
constant curvature spaces.  The integrand in Proposition \ref{cons}
now involves a sum over $I,I'$.  The general term inside the bracket
is  $c_{|I|}\epsilon(I,I')/[|I|!|I'|!]$ times
$$\sum_{\mu\in\Sigma_n}(\sgn\
 \mu)\frac{2^{|I|/2}\det(\frac{1}{2}(\Vert\circ df_x + \id))_I
\det(\frac{1}{2}(\Vert\circ df_x -
 \id))_{I'}\tilde\rho'(d(x,f(x))/\sqrt{2})}{[\cosh(d(x,f(x))/2)]^{|I|-1}
 \cdot [\sinh(d(x,f(x))/2)]^{|I'|-1}}$$
for negative curvature $-1$.  For constant
 curvature $1$,
 cosh, sinh are replaced by
cos, sin, and there is an extra factor of $(-1)^{|I|}$ due to
 $R_{ijij} = -1.$

\section{\large\rm\bf Large parameter behavior--topological methods}

In \cite[\S7]{MQ}, a one-parameter family of pullbacks of the
Mathai-Quillen form is
constructed which interpolates between the Hopf index formula (as
$t\to\infty$)
and the
Chern-Gauss-Bonnet theorem (at $t=0$).  
In this section, we show by topological arguments that the
corresponding 
large parameter behavior 
for Lefschetz theory is the Hopf fixed point/submanifold formula.
At the end, we indicate a geometric
refinement.  

For motivation, we sketch the argument in \cite{MQ}.  Let $s:M\to TM$
be a vector field on $M$ transverse to the zero section,
 and let ${\rm MQ}$ be the Mathai-Quillen form on
$TM$ for the Levi-Civita connection for a fixed Riemannian metric.
For  the zero
section $0$ of $TM$, we have $\Pf(\Omega) = 0^*{\rm MQ}$.
  For any $t\geq 0$, $\int_M
(ts)^*{\rm MQ}$ is independent of $t$,
 since the integrands are cohomologous.  As $t\to
\infty$, the integrand decays exponentially away from the fixed point
set, and the contribution from a fixed point $p$ becomes
$\pm\int_{T_pM} {\rm MQ} =1.$  Identifying the sign with the index
ind$(p)$ of
$s$ at $p$ gives 
\begin{equation}\label{pfhopf}
\int_M \Pf(\Omega) = 
\sum_{p, f(p) =p} {\rm ind}(p)\end{equation}
  These expressions equal $\chi(M)$ by
either Chern-Gauss-Bonnet or the Hopf index formula, but the
derivation of (\ref{pfhopf}) is new.

This argument still has content at the topological level:
if we replace MQ by any representative $\Phi$ of
the Thom class of $TM$,  we obtain at $t=0$ the well known expression
$\int_M 0^*\Phi$ for $\chi(M)$, and as $t\to\infty$ we still obtain
the
fixed point sum.  Thus Mathai and Quillen have produced a proof of the
Hopf index formula.  Note that the choice of section $s$ is
irrelevant.  In particular, the integral $\int_M s^*{\rm MQ}$ at $t=1$
has no particular significance.

To apply this method to the basic formula $L(f) = \int_M
(\id,f)^*{\rm MQ}$, we wish to replace $f$ by $tf.$  Since a tubular
neighborhood of the diagonal in $M\times M$ is diffeomorphic to $TM$,
 we can consider $f$ as a section of $TM$ whenever the graph of $f$
lies in the tube.  In this case,
we can define $tf$ as above. 
Moreover, $(\id,f)^*{\rm MQ}$ vanishes whenever the
graph lies outside the tube, so the action of $t$ might as well be
trivial on this set.  
This scaling of $f$ is done in detail below; for technical reasons we
deform the graph of $f$ in directions normal to the diagonal rather
than vertically.   Note that the
case $t=1$ now has particular significance: it is Theorem
\ref{theoremone}. 
In summary, we can think of $f$ as  a
section of $TM$ with possible blow up on part (or even all) of $M$.  

From this point of view, we can derive
the Lefschetz
formula for submanifolds of fixed points (Theorem 3.1) as
$t\to\infty.$
This proof is an elementary version of the
geometric stationary phase proof in \cite[Ch.~4]{G}.
To be
completely honest, there is a little geometry in the proof, but no
more
than in the usual proof of the tubular neighborhood theorem; even this
geometry can be eliminated by a jet bundle argument. 

In this section, $\nu_X^Y$ denotes the normal bundle of $X$ in $Y$.
\bigskip

To state the theorem, let $f:M\to M$ be a smooth map of a closed
oriented
$m$-manifold $M$,  and assume that the fixed
point set of $f$ consists of the disjoint union of smooth submanifolds
$N_j$ of dimension $n_j$.
Let $N$ be one such component, and let $\nu$ be the
quotient bundle $\nu=TM/TN$ over $N$. Since $df$ preserves the
subbundle $TN$, it induces a map $df_\nu$ on $\nu$.

We assume the non-degeneracy condition $\det(\id-df_\nu)\neq 0$ (also
known as clean intersection), i.e.~$f$ leaves infinitesimally fixed
only directions tangent
to $N$.

If we put $df_n$,  $n\in N$, in Jordan canonical form,  $TN$
will be the span
of eigenvectors with eigenvalues $1$, and $\nu$ is isomorphic to the
span of 
the generalized eigenvectors for the remaining eigenvalues. This
induces
a natural splitting of $TM\big |_N \simeq TN\oplus\nu.$ A choice of
Riemannian metric on $M$ gives an identification of $\nu$
with $\nu_N^M$, 
and $df_\nu$ with a map on $\nu_N^M$.

\begin{theorem}\label{theoremapp} Let $f:M\to M$ be a smooth
non-degenerate
map of a closed oriented $m$-manifold $M$,
whose fixed point set consists of the disjoint
union of  submanifolds $N_1, N_2, ..., N_r$.
Then
$$L(f)=
\sum_{j=1}^r\sgn(\det(\id -df_\nu))\chi(N_j).$$
\end{theorem}

To begin the proof, we may assume that the fixed point set of $f$
consists of a
single submanifold $N$ of dimension $n$.
Let $\Delta^\epsilon_M$ be an $\epsilon$-neighborhood of $\Delta_M$,
the
diagonal of $M$ in $M\times M$.   Choose $\epsilon >0$ small
enough so there exists a unique minimal geodesic from $x$ to $y$, for
all $(x,y)\in\Delta^\epsilon_M$. 

We construct a family of 
diffeomorphisms  
$F_t:M\times M\to M\times M$, for $t>0$, with $F_1 = \id$,  which
pushes out fibers of
$\nu_{\Delta_M}^{M\times M}$, while fixing $\Delta_M$ and $M\times
M-\Delta^\epsilon_M$.
Let $(x,y)\in\Delta^\epsilon_M$ and consider the geodesic
$\gamma$
in $M\times M$ from $(\bar{x},\bar{x})\in\Delta_M$ to $(x,y)$, 
where $\bar{x}$ is the midpoint of the geodesic $\alpha$ between $x$
and $y$ in $M$.
Setting 
$\dot\alpha(\bar x) = v= v(x,y)\in T_{\bar{x}}M$, we have
$\dot\gamma(\bar x,\bar x)=(-v,v).$
For $v\neq 0$, define a diffeomorphism
$\lambda_v(t):[0,\infty)\rightarrow [0,{\epsilon\over |v|})$  with
$\lambda_v(1)=1$, 
which is smooth in $v$,
and set $\lambda_0(t)=0$. Define
$F_t:M\times M\rightarrow M\times M$ by:
$$F_t(x,y)=\Bigg\lbrace\matrix{ (x,y),&
(x,y)\not\in\Delta^\epsilon_M,\cr
\exp_{(\bar{x},\bar{x})}(\lambda_{v(x,y)}(t)\cdot\exp^{-1}_{(\bar{x},\bar{x})}
(x,y)),& (x,y)\in\Delta^\epsilon_M.\cr}$$
$F_t$ is the desired map.  As in \S\ref{firstsection}, we have
$$L(f)=I(\Delta ,\Gamma)=
(-1)^{\dim\ M}I(\Gamma ,\Delta) 
=(-1)^{\dim\ M}\int_\Gamma\eta_{\Delta_M}^{M\times M},$$
where $\eta_{\Delta_M}^{M\times M}$
is the Poincar\'e dual of
$\Delta_M$ in $M\times M$. 
Since $F_t$ is homotopic to the identity we have
\begin{eqnarray}\label{starstar}
(-1)^{\dim\ M}L(f)
&=&\int_\Gamma\eta_{\Delta_M}^{M\times M}
= \lim_{t\rightarrow\infty}\int_\Gamma F^*_t\eta_{\Delta_M}^{M\times
M}
= \lim_{t\rightarrow\infty}\int_{(\id,f)(\Delta_M)} F_t^*
\eta_{\Delta_M}^{M\times M}\nonumber\\
&=&\lim_{t\rightarrow\infty}\int_{(\id,f)(\Delta^\delta_N)} F_t^*
\eta_{\Delta_M}^{M\times M}
= \lim_{t\rightarrow\infty}\int_{F_t\circ (\id,f)(\Delta^\delta_N)}
\eta_{\Delta_M}^{M\times M},\end{eqnarray}
where $\Delta^\delta_N$ is a $\delta$-neighborhood of
$\Delta_N$ in $\Delta_M$, for $\delta$ small enough.
This uses
$$\lim_{t\rightarrow\infty}\int_{(\id,f)(\Delta_M\setminus
\Delta^\delta_N)} F_t^*\eta_{\Delta_M}^{M\times M}=0,$$
as $F_t^*\eta_{\Delta_M}^{M\times M}$ decays uniformly
as $t\rightarrow\infty$ on $\Delta_M \setminus \Delta_N^\delta$,
since $d(x,f(x))$ and hence $|v|$ has positive minimum on $M$ minus
a $\delta$-neighborhood of $N$.

Let $\pi:\Delta_N^\delta\to\Delta_N$ be the projection given by
the identification of $\Delta_N^\delta$ with a
$\delta$-neighborhood of the zero section of 
$\nu_{\Delta_N}^{\Delta_M}$. 

The next lemma uses the non-degeneracy hypothesis.
\begin{lemma} If $df_\nu-\id$ is invertible at $n\in N$, then
$T_{(n,n)}\Big[F_t\circ(\id,f)(\pi^{-1}(n,n))\Big]\cap 
T_{(n,n)}\Delta_N=\{0\}$.\end{lemma}

\noindent {\sc Proof:}  If
$T_{(n,n)}\Big[F_t\circ(\id,f)(\pi^{-1}(n,n))\Big]
\cap T_{(n,n)}\Delta_N\neq{0}$,
there exists $0\neq
(q,q)\in T_{(n,n)}\pi^{-1}(n,n)$ with $q\perp N$
such that
$dF_t(q, df_\nu q) = dF_t\circ (\id,df)(q,q)\in T_{(n,n)}\Delta_N\ .$

We split  $(q,df_\nu q)\in T_{(n,n)}M\times M$ into its components in 
$T\Delta_M$ and in $\nu_{\Delta_M}^{M\times M}.$ 
 Since $dF_t$ leaves vectors in $T\Delta_M$ unchanged and
stretches vectors in the normal bundle by a $\lambda$ factor, we get
\begin{eqnarray*}
dF_t(q,df_\nu q)&=&dF_t
\bigg[\Big({q+df_\nu q\over 2}, {q+df_\nu q\over 2}\Big)
+\Big({q-df_\nu q\over 2}, {-q+df_\nu q\over 2}\Big)\bigg]\cr
&=&\Big({q+df_\nu q\over 2}, {q+df_\nu q\over 2}\Big)
+\lambda(t)\Big({q-df_\nu q\over 2}, {-q+df_\nu q\over
2}\Big)\cr
&=&\biggl({1+\lambda(t)\over 2}q+{1-\lambda(t)\over 2}df_\nu q,
{1-\lambda(t)\over 2}q+{1+\lambda(t)\over 2}df_\nu q\biggr),\cr
\end{eqnarray*}	
for $\lambda(t)=\lambda_w(t)$ with 
$w=(q-df_\nu q)/ 2.$ Note that by hypothesis, $w\neq 0$ and so
$\lambda(t) \neq 0.$

 We have
$dF_t\circ (\id,df)(q,q)=(v,v)$
for some $v\in T_nN$,
so
$$\Big({1+\lambda(t)\over 2}q+{1-\lambda(t)\over 2}df_\nu q\Big) -
\Big({1-\lambda(t)\over 2}q+{1+\lambda(t)\over 2}df_\nu
q\Big)=v-v=0.$$
This implies $\lambda(t)(\id-df_\nu)q=0.$
Since $q\neq 0,\  \lambda (t)\neq 0$, this contradicts 
that $\id-df_\nu$ is invertible.  \hfill$\Box$

\bigskip 

Define $E_{n,t}\subset
T_{(n,n)}(M\times M)$
by
$$E_{n,t} = T_{(n,n)}\Big[F_t\circ(\id,f)(\pi^{-1}(n,n))\Big],$$
 and note the decomposition
$$T(M\times M)\Big|_{\Delta_M}\simeq 
T{\Delta_M}\oplus\nu_{\Delta_M}^{M\times M}\simeq 
T\nu_{\Delta_M}^{M\times M}.$$
Let
$${\tilde \pi}:T(M\times M)\Big|_{\Delta_M}\rightarrow
\nu_{\Delta_M}^{M\times M}$$
be the projection to $\nu_{\Delta_M}^{M\times M}$.
By Lemma 2.1, ${\tilde\pi}$ has no kernel on $E_{n,t}$ and hence
is an isomorphism of $E_{n,t}$ to a vector subspace 
$H_{n,t}\subset\nu_{\Delta_M}^{M\times M}$. Let
$$\beta_{n,t}: E_{n,t}\rightarrow F_t\circ (\id,f)(\pi^{-1}(n,n))$$
be the diffeomorphism given by the exponential
map. Actually,
$\beta_{n,t}$ is a diffeomorphism on a neighborhood of $0$ in
$E_{n,t}$, whose radius goes to infinity 
as $t\rightarrow\infty.$

Thus, ${\tilde\pi}\circ\beta^{-1}_{n,t}: F_t\circ
(\id,f)(\pi^{-1}(n,n))
\rightarrow {\tilde H}_{n,t}\subset  H_{n,t}$ is a diffeomorphism onto
its
image ${\tilde H}_{n,t}$, where ${\tilde H}_{n,t}$ is an arbitrarily
large ball in $H_{n,t}$, for large $t$. Then
\begin{eqnarray*}(-1)^{\dim\ M}L(f)
&=&\lim_{t\rightarrow\infty}\int_{F_t\circ (\id,f)(\Delta^\delta_N)}
\eta_{\Delta_M}^{M\times M}\cr
&=&\lim_{t\to\infty}[{\rm deg}({\tilde
\pi}\circ\beta_t^{-1})^{-1}]^{-1}
\int_{({\tilde
\pi}\circ\beta_t^{-1})F_t\circ(\id,f)(\Delta^\delta_N)}
(({\tilde \pi}\circ\beta_t^{-1})^{-1})^*\eta_{\Delta_M}^{M\times
M}\cr
&=&\lim_{t\rightarrow\infty}{\rm deg}({\tilde
\pi})
\int_{\cup_n {\tilde H}_{n,t}} 
(({\tilde \pi}\circ\beta_t^{-1})^{-1})^*
\eta_{\Delta_M}^{M\times M},\cr
\end{eqnarray*}
where
$\beta_t: \cup_n E_{n,t}\rightarrow F_t\circ
(\id,f)(\Delta^\delta_N)$
is given by $\beta_{n,t}$ on each $E_{n,t}.$
This uses ${\rm deg}(\tilde\pi \circ \beta_t^{-1})^{-1} = ({\rm deg}
\ \tilde \pi)^{-1} = {\rm deg}\ \tilde \pi,$ as $\beta_t$ is an
orientation preserving diffeomorhism and $\tilde \pi$ is an
isomorphism.

Let $H_t=\cup_n H_{n,t}$ be the subbundle of $\nu_{\Delta_M}^{M\times
M}$  over $\Delta_N$
with
fiber $H_{n,t}$ over $(n,n)$. We obtain
$$(-1)^{\dim\ M}L(f)=\lim_{t\rightarrow\infty}
{\rm deg}({\tilde \pi})
\int_{H_t}(({\tilde \pi}\circ\beta_t^{-1})^{-1})^* 
\eta_{\Delta_M}^{M\times M}.$$

Since the $E_{n,t}$ are getting more ``vertical'' as
$t\rightarrow\infty$,
${\tilde\pi}\circ\beta_t^{-1}\rightarrow \pm \id$ and
$H_{n,t}\rightarrow
H_{n,\infty}$, where $H_{n,\infty}$ is the vector subspace of
$(\nu_{\Delta_M}^{M\times M})_{(n,n)}$
spanned by the projection of vectors in $H_{n,t}$ into
$\nu_{\Delta_M}^{M\times M}$, for any $t$.
 Set $H_\infty =\cup_n
H_{n,\infty}$
 with projection map
$ p:H_\infty\rightarrow\Delta_N.$  Then
$H_\infty$ is also a
subbundle
of $\nu_{\Delta_M}^{M\times M}\rightarrow \Delta_N$, and 
\begin{eqnarray}\label{appone}
(-1)^{\dim\ M}L(f)&=&\lim_{t\rightarrow\infty}
{\rm deg}({\tilde \pi})
\int_{H_t}(({\tilde\pi}\circ\beta_t^{-1})^{-1})^*
\eta_{\Delta M}^{M\times M}\nonumber\\
&=&{\rm deg}({\tilde \pi}) 
\int_{H_\infty}\eta_{\Delta_M}^{M\times M}.\end{eqnarray}

By Theorem \ref{theoremone}, $\eta_{\Delta_M}^{M\times M}$ and
$\Phi(  \nu_{\Delta_M}^{M\times M})$, the
Thom class of  $  \nu_{\Delta_M}^{M\times M}$
considered as a form on $M\times M$, can be
represented by the same form.  Since we do not distinguish between the
integral of a cohomology class and the integral of a representative
form, we have
\begin{eqnarray}\label{apptwo}
\int_{H_\infty}\eta_{\Delta_M}^{M\times M}
&=&\int_{H_\infty}\Phi(\nu_{\Delta_M}^{M\times M})
=\int_{H_\infty}\Phi(\nu_{\Delta_M}^{M\times M})\wedge 
p^*1\nonumber\\
&=&\int_{\Delta_N} p_*\Phi(\nu_{\Delta_M}^{M\times M})\wedge
1
=\int_{\Delta_N} p_*\Phi(\nu_{\Delta_M}^{M\times M}),
\end{eqnarray}
where the push forward formula for integration over the fiber
\cite[Prop.~6.15]{BT} is used
between the third and fourth terms.

Let $H_\infty^{\perp}$ be
the orthogonal (or any) complement of $H_\infty$ in $\nu^{M\times
M}_{\Delta_M}$.
By \cite[Prop.~6.19]{BT}, we have 
$\Phi(\nu^{M\times
M}_{\Delta_M})=\Phi(H_\infty)\wedge\Phi(H_\infty^\perp).$
It is easy to check that
\[
 p_*\Phi(\nu^{M\times M}_{\Delta_M})=
 p_*(\Phi(H_\infty)\wedge\Phi(H_\infty^\perp))
= p_*(\Phi(H_\infty))\wedge\Phi(H_\infty^\perp),\]
since $\Phi(H_\infty^\perp)$ vanishes in $H_\infty$ directions. Thus
 (\ref{apptwo}) becomes
\begin{equation}\label{appthree}
\int_{\Delta_N}p_*\Phi(\nu_{\Delta_M}^{M\times
M})
=\int_{\Delta_N}
p_*(\Phi(H_\infty))\wedge\Phi(H_\infty^\perp)
=\int_{\Delta_N}\Phi(H_\infty^\perp),
\end{equation}
as 
$p_*\Phi(H_\infty)=1$ since $\Phi(H_\infty)$ integrates to one in each
fiber.

We claim that $$\nu_{\Delta_M}^{M\times M}\Bigl|_{\Delta N}
\simeq\nu_{\Delta_N}^{N\times N}\oplus
\nu_{\Delta_N}^{\Delta_M}.$$
Indeed, the metric on $M$ is chosen
so that 
\begin{equation}\label{appfour}
TM\big|_N=TN\oplus\nu_{N}^{M}.\end{equation}
$\nu_{\Delta_N}^{N\times N}$ is isomorphic to 
$TN$ by the map $(v,-v)\mapsto v.$  Similarly,
$\nu_{\Delta_M}^{M\times M}\simeq T\Delta_M\simeq TM.$  Finally, we
trivially have
$\nu_{\Delta_N}^{\Delta_M}\simeq \nu_{N}^{M}.$  Plugging these terms
into (\ref{appfour}) gives the claim.

Thus we have the bundle isomorphisms
$\nu_{\Delta_N}^{\Delta_M}\simeq E_t\simeq H_t\simeq H_\infty,$ and by
the claim we have
$H_\infty^\perp\simeq\nu^{N\times N}_{\Delta_N}.$   By (\ref{apptwo}),
(\ref{appthree}), we have
\begin{eqnarray}\label{appfive}\int_{H_\infty}\eta_{\Delta_M}^{M\times
M}&=&
\int_{\Delta_N}\Phi(H_\infty^\perp)
=\int_{\Delta_N}\Phi(\nu_{\Delta_N}^{N\times N})
=\int_{\Delta_N}\eta_{\Delta_N}^{N\times N}\nonumber\\
&=&I(\Delta_N, \Delta_N)
=\chi(N),\end{eqnarray}
where the self-intersection number of $\Delta_N$ appears as
in \S2.1.
Combining (\ref{appone}),  (\ref{appfive}) gives the
Lefschetz formula up to sign:
\begin{equation}\label{appsix}L(f) = (-1)^{ {\rm dim}\ M}{\rm
deg}(\tilde \pi)\chi(N).\end{equation}

To compute the degree of  ${\tilde\pi}: E_{n,t}\rightarrow 
H_{n,t}\simeq  \nu_{\Delta_N}^{\Delta M}$, we
pick $\theta$ and $\alpha$,
 positively oriented bases for $E_{n,t}\subset T_{(n,n)}\Gamma$
and $H_{n,t}\simeq (\nu_{\Delta_N}^{\Delta M})_{(n,n)}$
respectively,
and compute the sign of the determinant of the matrix of ${\tilde\pi}$
with respect to $\theta$,
$\alpha$. 

There exists a positively oriented basis for $T_nM$, 
$(v_1,... ,v_n, w_{n+1},... ,w_m)$, with $v_i\perp w_j$, 
such that
$v_1,... ,v_n\in T_nN,\ df_n v=v$ and
$w_1,... ,w_{m-n}\in \nu_N^M,\ df_n w=df_{\nu_n} w$.
A positively oriented basis for $T_{(n,n)}\Gamma$ is then 
$$\lbrace (v_1,v_1),... ,(v_n,v_n), (w_1, df_\nu w_1), ...
,... ,(w_m,df_\nu w_{m-n})\rbrace,$$
and a positively oriented basis for $E_{n,t}$ is
$$\theta =\lbrace (w_1,df_\nu w_1),... 
,... ,(w_m,df_\nu w_{m-n})\rbrace,$$
since $E_{n,t}\simeq (\nu_{\Delta_N}^{\Delta_M})_{(n,n)}$.
A positively oriented basis for 
$H_{n,t}\simeq (\nu_{\Delta_N}^{\Delta_M})_{(n,n)}$ is
$$\alpha =\lbrace (-w_1,w_1),... ,(-w_{m-n},w_{m-n})\rbrace.$$
As in Lemma 2.1, the vectors in $\theta$ decompose into
$$(w_i, df_\nu w_i)=\Big({w_i+df_\nu w_i\over 2}, {w_i+df_\nu w_i\over
2}\Big)
+\Big({w_i-df_\nu w_i\over 2}, {-w_i+df_\nu w_i\over 2}\Big).$$
Hence
\begin{eqnarray*} {\rm deg}\ {\tilde\pi}
&=&\sgn\ \det\bigg\lbrace(-w_i,w_i)
\mapsto \bigg({w_i-df_\nu w_i\over 2}, 
{-w_i+df_\nu w_i\over 2}\bigg)\bigg\rbrace\cr
&=&\sgn\ \det\bigg\lbrace(-w_i,w_i)\mapsto
(df_\nu-\id)(-w_i,w_i)\bigg\rbrace\cr
&=&\sgn\ \det (df_\nu-\id).\end{eqnarray*}
Since the right hand side of (\ref{appsix}) vanishes if dim $N$ is
odd, we assume dim $N$ is even.  (\ref{appsix}) becomes
\begin{eqnarray*} L(f) &=& (-1)^{{\rm dim}\ M} \sgn
(\det(df_\nu-\id))\chi(N)\\
&=& (-1)^{{\rm dim}\ M} (-1)^{{\rm dim}\ M - {\rm dim} \ N}\sgn (\det
(\id - df_\nu))\chi(N)\\
&=& \sgn (\det
(\id - df_\nu))\chi(N),\end{eqnarray*}
which concludes the proof of Theorem 3.1.

If the graph of $f$ is transversal to the diagonal, 
the fixed point set reduces to a finite
number of isolated fixed points $n_1, n_2,... ,n_r$, and the Lefschetz
fixed point formula is easily recovered.
For let $n$ be one such isolated fixed point. Then $H_\infty$
reduces
to $H_{n,\infty}$, the fiber over $(n,n)$ in $\Delta_N$ and
$\int_{H_{n,\infty}}\Phi(\nu_{\Delta_M}^{M\times M})=1$.
$df_\nu$ is just $df_n$
and
${\rm deg}({\tilde\pi})=\sgn\ \det(
df_n-\id).$
So (\ref{appone}), (\ref{apptwo}) give the fixed point formula
$$L(f) = (-1)^{{\rm dim}\ M}\sum_{i=1}^r {\rm deg}(\tilde
\pi_i)\int_{H_{n_i,\infty}}\Phi(\nu_{\Delta_M}^{M\times M}) =
\sum_{i=1}^r {\rm sgn}\  \det(\id -df_{n_i}).$$

\bigskip

\noindent{\bf Remark:}  We sketch the corresponding geometric proof
of
the fixed submanifold formula based
on Theorem \ref{localtheorem}.
Assume that the metric on $M$ is a product
near a fixed point submanifold $N$. If the submanifold is given by
$\{
x^{k+1} = \ldots = x^n =0\}$ in local coordinates, then as
$t\to\infty$, the integrand for $L(f)$ concentrates on a tubular
neighborhood of the fixed point, and the only contribution to the
integrand comes from $I = \{1,\ldots,k\}$, since the curvature term
vanishes otherwise due to the product metric.  Converting back to
rectangular coordinates in the normal fiber as in the topological
proof eliminates the $\tilde\rho^\prime$ factor and introduces
a factor of sgn det$(df_\nu - {\rm Id})$. Since $f={\rm Id}$ in
submanifold
directions, $\det(
\frac{1}{2}(d(tf)+{\rm Id})_{I}^\mu )= 1$.  Thus the integral splits
into
the curvature integral over $N$, yielding $\chi(N)$, and a normal
integral, which gives sgn det$(d(tf)_\nu - {\rm Id})$.  In the
$t\to\infty$ limit, $d(tf)_\nu-{\rm Id}$ in the normal fiber goes to
the
identity map, so its determinant becomes one.  Plugging these terms
into the integrand in Theorem \ref{localtheorem} gives
the Lefschetz fixed submanifold formula.

\section{\large\rm\bf Small parameter behavior--the role of the cut locus}

For a section $s$ of $TM$, the behavior of $ts$ as $t\to 0$ is
trivial. In contrast, the function $tf$ mentioned in the beginning of
\S3 (or more precisely, the diffeomorphism $F_t$)
 becomes discontinuous at $t=0$ at those $x\in M$ for which
 $(x,f(x))$ is
in the boundary of the tubular neighborhood of the diagonal.  
If the graph lies entirely within the tube, then $tf$ is well defined,
$\lim_{t\to 0}tf = \id,$ and $L(f) = \chi(M)$.

Thus we expect the
difference between $L(f)$ and $\chi(M)$
to be concentrated on the intersection of the graph with the
tube boundary.  Since this intersection is quite bad in general, the
difference will be given by a current supported on the intersection.

The maximum amount of information is obtained when the tube is as
large as possible.  As explained below,
this occurs when the boundary of the vertical
fiber of the tube at $(x,x)$ is $\cutx$, the cut locus of $x$.  Recall that
on a closed manifold, a geodesic $\gamma(t)$ is the minimal length
curve joining $x=\gamma(0)$ to $\gamma(t)$  for $t\in [0,T]$ for some
finite time $T$.  The point $y = \gamma(T)$ is by definition in
$\cutx$.

A point $y\in \cutx$ is characterized by: either there
exists more than one minimal length geodesic from $x$ to $y$, or
$d\exp_x$ is singular at the preimage of $y$ \cite[Lemma 5.2]{CE}.
In particular,
if the graph of a smooth function $f:M\to M $ has the property that
$f(x)$ is never on $\cutx$, then there is a
unique minimal geodesic joining $x$ to $f(x)$.
Shrinking this geodesic gives a homotopy from $f$ to the identity, and
so the Lefschetz number satisfies $L(f) =
\chi(M).$    Thus the difference $L(f) -\chi(M)$ is controlled by the
{\it cut locus of $f$ in $M$}
\begin{equation}\label{cutf} \cutf =
\{x:f(x) \in \cutx\}.\end{equation}

  In the first subsection, we will make this
geometric statement more precise by finding a singular current supported on the
cut locus of $f$ whose singular part evaluated at the function $1$ gives $L(f)
-\chi(M)$.  The main idea is to define the function $tf$ and to let
$t\to 0.$ 
In the second subsection, we assume that $\cutf$ is finite.
Under a transversality
condition, 
the number of points in $\cutf$ can be estimated from below.
In fact, for
all but very special metrics, the transversality condition
implies that $\cutf$ is infinite for diffeomorphisms with
 $L(f) \neq \chi(M).$

\subsection{\large\rm\bf A current on the cut locus}

We first construct the largest (topological) tubular neighborhood of 
the diagonal $\Delta$  in $M\times M$.  A tubular
neighborhood is
given by points of the form $(\exp_{\bar x}v,\exp_{\bar x} (-v))$,
where $v\in T_{\bar x}M, \bar x \in M$, and $|v|$ is small.
For $x,y\in M$, we say that $y$ is
inside $\cutx$ if there is a unique minimal geodesic from $x$ to $y$.
Let $N_{\bar x} = \{\exp_{\bar x}v: \expx v\ {\rm
is\  inside}\ {\cal C}_{\expx (-v)}\}.$   
Define $T\subset M\times M$ by $T = \{\evev:\expx v \in N_{\bar x},
\bar
x\in M\}.$  We call $T$ the {\it cut locus tubular neighborhood.}

\begin{lemma}\label{firstlemma_a}
(i)  $T$ is a topological tubular neighborhood of the
diagonal.  

(ii) $(x,y) \in T$ iff $y$ is inside $\cutx.$

(iii) The vertical fiber $T\cap (\{x\}\times M)$ of $T$ at $x$ equals
$\{x\} \times (M\setminus \cutx).$
\end{lemma}

\noindent{\sc Proof:}  We prove (ii)  first.  The
forward implication is from the definition
of $T$.  Conversely, if $y$ is inside $\cutx$, then there is a unique
minimal geodesic from $x$ to $y$.  Then
$(x,y) = \evev$, where $\bar x$ is the midpoint of the geodesic.  Thus
$x, y\in N_{\bar x}$, so $(x,y) \in T.$

For (i), the standard argument that the interior of
the cut locus is a
topological sphere immediately extends to show that
$\exp_x^{-1}(N_x)\subset T_xM$ is
the
interior of a topological sphere.  This argument in turn extends to
show that the radius of this sphere is a continuous function on the
unit tangent sphere, which implies that $\exp^{-1}(T)$ is a topological disk
bundle.

To finish the proof, we must show that $\exp$ on $M\times M$ is
injective on $\{(v,-v)\}\subset \coprod_x(\exp_x^{-1}(N_x)\times
\exp_x^{-1}(N_x)).$   If not, there exists $\bar x,\bar y,v,w$ with
$\alpha = \expx v = \exp_{\bar y} w $ and $\beta = \expx (-v) =
\exp_{\bar y} (-w).$  By the definition of $N_{\bar x},N_{\bar y}$,
this gives two minimal geodesics from $\alpha$ to $\beta$, a
contradiction. 

For (iii), $(\exp_{\bar x}v,\exp_{\bar x} (-v))$ is in the vertical
fiber over $(\exp_{\bar x}v,\exp_{\bar x} v)$, and at $\partial T$,
$\exp_{\bar x} (-v) \in {\cal C}_{\exp_{\bar x}v}.$
\hfill$\Box$\break

\bigskip


We now define $tf$. Fix a diffeomorphism $\mu:[0,1)\to[0,\infty)$ with
$\mu(0) = 0,\mu(1) = 1$ and such that the  derivative of
$\mu^{-1}$ grows at most polynomially.  For the moment, let $T$ denote
any smooth tubular neighborhood of $\Delta$ given by geodesics of the
form
$(x,y) = (\gamma(t),\gamma(-t)).$ (The range of $t$ is a smooth
function on the unit tangent bundle.)
 For such $(x,y)$, set
$$d_{x,y} = \min \{t:(\gamma(t),\gamma(-t)) \in \partial
T\}.$$    For $x\in M$ and $t\in [0,\infty)$, define $t_x:M\to M$ by
$$t_x(y) = \left\{ \begin{array}{ll}
  \exp_x[\mu^{-1}(\mu(d_{x,y}^{-1}|v|)t)\dxy
\frac{v}{|v|}],& (x,y)\in T, \ y = \exp_x v, \ y \neq x,\\
y,& (x,y) \not\in T,\\
x,& y=x.  \end{array} \right.$$
Thus for $x,y$ close, $t_x$ pushes $y$ towards $\partial T$ as
  $t\to\infty$
 along their minimal geodesic,
but fixes $y$ if it is far from $x$, as measured by $T$.

For $f:M\to M$, define $tf:M\to M$ by
$$(tf)(x) = t_x(f(x)).$$
The maps $tf$ are smooth for $t>0.$  
Note that $(1f)(x) = \exp_x v = f(x)$ if $(x,f(x))\in T$ and $(1f)(x)
=
f(x) $ otherwise, so $1f = f.$  Similarly, we have $(0f)(x) = x$ if
$(x,f(x))\in T$, and $(0f)(x) = f(x)$ otherwise.  
Thus $0f$ is discontinuous on $\{x:(x,f(x))\in \partial T\}.$  We
remark that $(\id, tf)$ can be used in place of $F_t$ to prove the
Lefschetz fixed submanifold formula as $t\to\infty.$

We now examine the $t\to 0$ limit of pullbacks of Mathai-Quillen
forms.
Fix $\epsilon$, and let ${\rm MQ}_\Delta = {\rm MQ}_{\Delta_\epsilon}$
be the Mathai-Quillen form on the $\epsilon$-neighborhood of the
diagonal. Then
\begin{eqnarray*} L(f) &=& \int_\Delta ({\rm Id},f)^* \MQ_\Delta\\
&=& \lim_{t\to 0} \int_\Delta ({\rm Id},tf)^* \MQ_\Delta\\
&=& \lim_{t\to 0} \int_{\{(x,x):(x,f(x))\in T\}} ({\rm Id},tf)^*
\MQ_\Delta +
\lim_{t\to 0} \int_{\{(x,x):(x,f(x))\not\in T\}}  ({\rm Id},tf)^*
\MQ_\Delta.\end{eqnarray*}
 As usual $({\rm Id},tf)^* (\MQ_\Delta)_{(x,x)} =
(\MQ_\Delta)_{(x,f(x))} \circ ({\rm Id},tf)_* = 0$ if $(x,f(x))\not\in
T,$  so
$$L(f) = \lim_{t\to 0} \int_{\{(x,x):(x,f(x))\in T\}} ({\rm Id},tf)^*
\MQ_\Delta.$$ 
Note that $d(x,y)/d_{x,y} <2|v|/|v| = 2$ for $(x,y) = (\exp_{\bar x}v,
\exp_{\bar x}(-v)).$
Fix $\delta <1$, set 
$${\adel} = \{x:(x,f(x))\in T, \frac{d(x,f(x))}{d_{x,f(x)}} \leq
2\delta\},$$
and set ${\bdel} = \{x:(x,f(x)) \in T\} \setminus {\adel}.$  It is
easy to
check that $({\rm Id},tf)^*\MQ_\Delta \to ({\rm Id},{\rm
Id})^*\MQ_\Delta = {\rm Pf}(\Omega)$ as $t\to 0$ uniformly on the
compact set
${\adel}.$  Thus
\[L(f) = \int_{\adel} {\rm Pf}(\Omega) + \lim_{t\to 0} \int_{\bdel}
i^*\itfmq,\]
where $i:M\to \Delta$ is the inclusion; we will omit this map from
here on.
Since ${\rm Pf}(\Omega)$ is smooth on $M$ and ${\adel}$ exhausts the
open set 
$\{x:(x,f(x))\in T\}$ as $\delta\to 1$, we get
\begin{equation}\label{star}
L(f) = \int_{\{x:(x,f(x))\in T\}} \euler + \lim_{\delta\to
1}\lim_{t\to 0} \int_{\bdel}\itfmq.\end{equation}

This construction extends to the topological tubular neighborhood $T$ of
Lemma \ref{firstlemma_a}.  
Fix $\epsilon >0$ and pick a smooth disk
bundle $D^\epsilon\subset \nu_\Delta$ such that $T^\epsilon =
\exp D^\epsilon$ is
inside $T$ and is within $\epsilon$ of filling $T$ -- i.e.~for all
$(v,-v)\in \partial D^\epsilon$, we have $d_{M\times M}( (\exp v,
\exp(-v)), (\exp(tv),\exp(-tv))) <\epsilon$, where
$t$ is the smallest positive number such that $(\exp
tv, \exp(-tv))\in \partial T.$  
To define the Mathai-Quillen form on $T^\epsilon$, we have to choose a
diffeomorphism 
$\alpha^\epsilon:\nu_\Delta\to D^\epsilon$ and pull back the
Mathai-Quillen form $\MQ_\nu$ from $\nu_\Delta.$
 As $\epsilon\to 0$, $D^\epsilon$ fills out a continuous disk bundle
in $\nu_\Delta,$ and we demand that for all $R>0, $ there exists
$\epsilon_0 = \epsilon_0(R)$ such that
$\alpha^\epsilon(B_R(\nu_\Delta))$ is constant for all $\epsilon
<\epsilon_0$, where $B_R(\nu_\Delta)$  is the $R$-ball around the zero
section in $\nu_\Delta.$  For this choice of $\alpha^\epsilon$, it is
immediate that 
$$\MQ^0_\Delta(v_1,\ldots,v_n) \equiv \lim_{\epsilon\to
0}[((\alpha^\epsilon)^{-1})^*\MQ_\nu(v_1,\ldots,v_n)]$$
exists and is smooth, and that $\MQ_\Delta^\epsilon \equiv
((\alpha^\epsilon)^{-1})^*\MQ_\nu
\to \MQ^0_\Delta$ pointwise.  In fact, since $\MQ_\nu$
 decays exponentially at infinity in $\nu_\Delta$, it is easy to check
that this convergence is uniform. This yields
$$L(f) = \lim_{\epsilon\to 0}\int_\Delta \itfmq^\epsilon = 
\int_\Delta \itfmq^0.$$

We can now repeat the argument leading to (\ref{star}), noting that
for the topological neighborhood $T$, $d_{x,y}$ is just continuous in
$x,y.$  By Lemma \ref{firstlemma_a}, we obtain
\begin{equation}\label{ezero} 
L(f) = \int_{M\setminus \cutf} \euler + \lim_{\delta\to 1}\lim_{t\to
0}
\int_{\bdel}\itfmq^0,\end{equation} 
where $\cutf$ is given by (\ref{cutf}).  Note that $\cutf$ is closed,
so the first integral exists.
By the Chern-Gauss-Bonnet theorem, we get
\begin{equation}\label{eone}
 L(f) = \chi(M) - \int_{\cutf} \euler + \lim_{\delta\to 1}\lim_{t\to
 0} \int_M
\chi_{_{\bdel}}\cdot \itfmq^0.\end{equation}

We now define zero currents $L^{tf}, E$,
on $M$ via their action on $g\in
C^\infty(M)\simeq C^\infty(\Delta)$:  
\begin{eqnarray*} L^{tf}(g) &=& \int_M g\cdot\itfmq^0,\\
E(g) &=& \int_M g\cdot\euler.\end{eqnarray*}
We also set
\[
{\cal C}^f(g) = -\int_{\cutf }g\cdot\euler + 
\lim_{\delta\to 1}\lim_{t\to 0} \int_M
g\cdot\chi_{_{\bdel}}\cdot\itfmq^0,\]
whenever the right hand side exists.
  We define the limit of currents by pointwise convergence:
  $\lim_{t\to 0} L^{tf} = L^0$ if
  $\lim_{t\to 0} L^{tf}(g) = L^0(g)$ for all smooth $g$.

\begin{lemma} As a current, $(\lim_{t\to 0}
 L^{tf}) - \cf$ exists and
equals $E$.  In particular, $\lim_{t\to 0} L^{tf}
(g)$ exists whenever
{\rm supp} $g\cap \cutf = \emptyset.$
\end{lemma}

\noindent{\sc Proof:}
We have
\[L^{tf}(g) =  \int_M g\cdot\itfmq^0
= \int_{\adel} g\cdot\itfmq^0 + 
\int_{\bdel}g\cdot\itfmq^0,\]
and so
\begin{eqnarray*} \lim_{t\to 0} L^{tf}(g) - \lim_{\delta\to 1}\lim_{t\to 0}
  \int_{\bdel}g\cdot\itfmq^0 &=&
 \int_{M\setminus\cutf}
g\cdot\euler\\ 
&=& E(g) - \int_{\cutf} g\cdot\euler,\end{eqnarray*}
as in (\ref{ezero}).  This 
gives the first statement.  For the second statement, if ${\rm supp}\ g\cap
\cutf = \emptyset$, then $g\cdot\chi_{D_\delta} =0$ for $\delta\approx
1$, and so ${\cal C}^f(g) =0$ for such $g$. \hfill$\Box$\break
\bigskip

In view of this lemma, we think of $L^0$ as a
singular current, with $\cf$  the singular part of 
$ L^{0}$ and $E$ the finite part.
Note that $L^{tf}(1) = L(f)$ for all $t.$  This gives:

\begin{theorem}\label{singcur}
  For every Riemannian metric on a closed manifold
$M$ and every smooth function $f:M\to M$,
there exists a canonical singular part ${\cal C}^f$ to 
$L^0 = \lim_{t\to 0} L^{tf}$, with supp ${\cal C}^f
\subset \cutf$. Moreover, we have
\[ L(f) = \chi(M) + {\cal C}^f(1).\]   \end{theorem}

\noindent{\bf Remarks:}

(1)  The previous discussion can be watered down to apply to the
degree of
$f$,  defined by $\degf = \int_M f^*\omega/\int_M\omega$ for
any top degree form $\omega.$   By its homotpy invariance,
the degree of $f$ is one if the
graph of $f$ never intersects $\cutx$, so we expect that a singular
0-current, with singular part supported on $\cutf$, computes
$\degf - 1.$
Taking $\omega$ to be the volume form
dvol of a Riemannian metric on $M$, we get
\begin{eqnarray*} {\rm vol}(M)\cdot\degf &=& \lim_{t\to 0} \int_M
(tf)^*{\rm dvol}\\ 
&=& \int_{M\setminus \cutf} {\rm dvol} + 
\lim_{\delta\to 1}\lim_{t\to 0} \int_M
\chi_{_{\bdel}} (tf)^*{\rm dvol}.\end{eqnarray*}
Setting 
\[{\cal D}^f(g) = 
-\int_{\cutf }g\cdot {\rm dvol} + 
\lim_{\delta\to 1}\lim_{t\to 0} \int_M
g\cdot\chi_{_{\bdel}}\cdot (tf)^*{\rm dvol},\]
whenever the right hand side exists,
we see that the support of ${\cal D}^f$ is contained in $\cutf$, and
that
\[\degf -1 = \frac{{\cal D}^f(1)}{{\rm vol}(M)}.\]

(2) In (1) and in the previous section, we have compared $f$ to (the
homotopy class of) the identity map.  We can also compare $f$ to a
constant map $c(x) = x_0.$  In this case the Lefschetz number (resp. 
degree) of $f$ is 1 (resp. 0) if the graph of $f$ misses ${\cal
C}_{x_0}.$  Again there are singular 0-currents, with singular part
supported on ${\cal
C}_{x_0}$, which measure $L(f) -1$ and $\degf.$  

(3)  Finally, we can compare $f$ to a fixed map $f_0.$  We obtain
\[L(f) - L(f_0) = {\cal E}^{f,f_0}(1),\]
where the singular 0-current ${\cal E}^{f,f_0}$ has singular part
supported on $\{x:f_0(x) \in
{\cal C}_{f(x)}\}.$  There is a similar result for degrees.  As a
simple well-known 
example, note that if $M= S^n$ and $f,f_0$ have different
Lefschetz numbers (equiv. different degrees), then there exists $x\in
S^n$ such that $f(x), f_0(x)$ are antipodal.

\bigskip

$\cf(1)$ can be identified in the simplest
case where $f$ is Lefschetz (i.e.~the graph $\Gamma$ of $f$ is
transverse to $\Delta$--a generic condition) and
$\cutf$ consists of isolated points $\{x_1,\ldots,x_n\}$
(i.e.~$\Gamma\cap \partial T = \{(x_i,f(x_i))\}$, which we will see is
a non-generic condition).
Since $T$ is homeomorphic to $TM$, and diffeomorphic away from
$\partial T$, we can
consider $f$ as a smooth vector field $V_f$ on $M$ with singularities at the
$x_i$.
At each fixed point $x$ of $f$, the local Lefschetz number $L_x(f)$
equals the Hopf
index ind${}_x(V_f)$ 
of $V_f$ \cite[p.~135]{GP}.  We modify $V_f$ by multiplying  the
vectors in a neighborhood of each $x_i$ by a smooth function which is
one on the boundary of the neighborhood and which vanishes to all
orders at $x_i.$   The modified vector
field $\vfp$ extends to a smooth vector field, also denoted $\vfp$,
on all of $M$ with zeros at the fixed points of $f$ and at the $x_i.$
We have
$$\chi(M) = \sum_{\{x:\vfp(x) =0\}} {\rm ind}_x(\vfp) =
\sum_{\{x:f(x) =x\}}L_x(f) + \sum_i {\rm ind}_{x_i}(\vfp)
 = L(f) + \sum_i {\rm ind}_{x_i}(\vfp).$$

\begin{proposition} Let $f:M\to M$ be a Lefschetz map with $\cutf$
consisting of isolated points.  Then
\[
\cf(1) = -\sum_{x_i\in \cutf} {\rm ind}_{x_i}(\vfp),\]
and in particular
\[L(f) = \chi(M) - \sum_{x_i\in \cutf} {\rm ind}_{x_i}(\vfp).\]
\end{proposition}

\noindent{\bf Remark:} This proposition is related to the proof
of the Hopf index formula in \cite{MQ}.  As in the beginning of \S3,
for a vector
field $s$ we have
$\chi(M) = \int_M (ts)^* {\rm MQ}_{TM}$ for all $t$.  At $t=0$ we recover
the Chern-Gauss-Bonnet formula, so
$\chi(M) = \lim_{t\to\infty} \int_M (ts)^* {\rm MQ}_{TM}.$
Let $B_\epsilon$ be the $\epsilon$-neighborhood of the zero set of
$s$.  Then by the uniform decay of $(ts)^\ast {\rm MQ}_{TM}$ off
$B_\epsilon$, 
$\chi(M) = \lim_{\epsilon\to
0}\lim_{t\to\infty}\int_{B_\epsilon}(ts)^* {\rm MQ}_{TM}.$  
Define a family of Euler
currents by
$$E_s^t(g) = \int_M g\cdot(ts)^* {\rm MQ}_{TM}.$$
Then $E_s^t(1) = \chi(M)$ and
$$\lim_{t\to \infty} E_s^t -
\lim_{\epsilon\to
0}\lim_{t\to\infty}\int_{B_\epsilon}(ts)^* {\rm MQ}_{TM} =0.$$
Thus the singular 0-current $\lim_{t\to\infty}E_s^t$ is supported on the
zero
set of $s$, and the Euler characteristic as a 0-current
localizes to 
the zero set.  If the zero set
consists of nondegenerate points, this singular part is given by
$\pm\delta$-functions at the zeros, and the Hopf index formula is
recovered.

\subsection{\large\rm\bf Isolated cut points}

Assume that (i) $\cutf$ consists of a finite set of points, and
(ii) the graph
$\Gamma$ of $f$ 
is transverse to $M\times \{f(x)\}$ for all $x\in \cutf.$  In this
case, we say that {\it f is transverse to the cut locus.}  
Under this assumption, we will show that
$|\cutf|$ can be bounded from below.

Condition (ii) is equivalent to $df_{x}$ being invertible, as $(v,
df_{x}v) \in T(M\times \{f(x)\})$ implies $df_xv=0.$  In particular, a
diffeomorphism of $M$ satisfies (ii).

For simplicity, write $\cutf = \{x\}.$
The transversality assumption
 implies that the differential of $f$ is invertible at $x$, so on
some $\epsilon-$neighborhood $B_\epsilon(x)$,
$\Gamma|_{B_\epsilon(x)}$ is a
graph over the neighborhood $U = \{x\}\times f(B_\epsilon(x))$ of
$(x,f(x))\in \{x\}\times M.$   Thus the projection $p_2:M\times M\to
M$ onto the second factor restricts to a diffeomorphism
$p_2:\Gamma|_{B_\epsilon(x)}\to U,$ which of course has degree $\pm
1.$
For ${\rm MQ}_\Delta$ the Mathai-Quillen form of
$\nu_{\Delta}^{M\times M}$, considered as a form on the cut locus
tubular neighborhood, we have
\begin{eqnarray} \label{one_a}
L(f) &=& \lim_{\epsilon\to 0} \lim_{t\to 0}
\int_{M\setminus B_\epsilon(x)}({\rm Id}, tf)^* \MQ_{\Delta} + 
\lim_{\epsilon\to 0} \lim_{t\to 0}
\int_{B_\epsilon(x)}({\rm Id}, tf)^* \MQ_{\Delta}\nonumber\\
&=& \lim_{\epsilon\to 0}\int_{M\setminus B_\epsilon(x)} \euler +
\lim_{\epsilon\to 0} \lim_{t\to 0}
\int_{B_\epsilon(x)}({\rm Id}, tf)^* \MQ_{\Delta}\\
&=& \chi(M)  + 
\lim_{\epsilon\to 0} \lim_{t\to 0}
\int_{B_\epsilon(x)}({\rm Id}, tf)^* \MQ_{\Delta}.\nonumber
\end{eqnarray}
Here we do not distinguish between integrals over $M$ and
integrals over the diagonal $\Delta\subset M\times M.$
To justify (\ref{one_a}), we need that $({\rm Id}, tf)^*
\MQ_{\Delta}$ converges uniformly to $\euler$ on $M\setminus
B_\epsilon(x)$ as $t\to 0$.

\begin{lemma} Let $\imath:\Delta \to M\times M$ be the inclusion.
Then
\[\lim_{t\to 0} \imath^*({\rm Id},tf)^* \MQ_\Delta = \euler\]
uniformly on $M\setminus B_\epsilon (x).$ \end{lemma}

\noindent {\sc Proof:} On $M\setminus B_\epsilon (x),$ we have
$(tf)(y) \to y$ uniformly as $t\to 0$.  Thus if $\gamma(s)$ is a short
curve with $\gamma(0) = y, \dot\gamma(0) = w$, then $(tf)(\gamma(s))
\to \gamma(s)$ uniformly as $t\to 0,$ and
\begin{eqnarray*} \lim_{t\to 0} (tf)_*(w) &=& \lim_{t\to 0} \lim_{s\to
0} \frac{(tf)(\gamma(s)) - (tf) (y)}{s} \\
&=& \lim_{s\to 0} \lim_{t\to
0} \frac{(tf)(\gamma(s)) - (tf) (y)}{s} \\
&=& \lim_{s\to 0}\frac{\gamma(s) - \gamma(0)}{s} = \dot\gamma(0)
=w.\end{eqnarray*} 
This shows that
\begin{eqnarray*} 
\lefteqn{[\itfmq]_{(y,y)}((v_1,w_1),\ldots,(v_n,w_n))= }\\
&&  (
\MQ_{\Delta})_{(y,f(y))}((v_1,(tf)_*w_1),\ldots,(v_n,(tf)_*w_n))
\end{eqnarray*}
converges uniformly in $y$ to
\[ ({\rm
MQ}_{\Delta})_{(y,y)}((v_1,w_1),\ldots,(v_n,w_n))\]
as $t\to 0.$  Since $\imath^*{\rm MQ}_\Delta = \euler$, the lemma
follows.
\hfill$\Box$

\bigskip
Let $\tmup$ denote $\{0\}\times TM\subset T(M\times
M)|_\Delta,$  
and let ${\rm MQ}_{\tmup}$ denote the Mathai-Quillen form of $\tmup$,
considered as a form supported on the cut locus neighborhood.  Thus
${\rm MQ}_{\tmup} = (\exp^{-1})^*\beta^*{\rm MQ}$, where MQ is the
Mathai-Quillen form on $\tmup$, 
$\exp$ is the exponential map from $\tmup$ to $M\times M$, and
$\beta$ is a 
homeomorphism from the neighborhood of zero in $\tmup$ with fiber 
$\exp_x^{-1}(M\setminus \cutx)$ 
to $TM$.  Here we have used Lemma \ref{firstlemma_a} (iii).
As in the last section, $\beta$ is a limit of
diffeomorphisms, and because of the decay of MQ we may treat $\beta$
as a diffeomorphism.  

Note that $p_2^* \MQ_{TM^\uparrow} = \MQ_\Delta$, since 
(i) $p_2^* \MQ_{TM^\uparrow} $ is closed and (ii) for a
fiber $F = \{(\exp_x (-v), \exp_x v): v\in N_x\}$ of the cut locus
tubular neighborhood, we have 
\begin{eqnarray}\label{yah}
\int _F p_2^*{\rm MQ}_{\tmup} &=& \int_{p_2F}{\rm MQ}_{\tmup}
=\int_{M\setminus \cutx} {\rm MQ}_{\tmup}\nonumber\\
&=& \int_{\beta\exp_x^{-1}(M\setminus \cutx)} {\rm MQ}
= \int_{T_xM^\uparrow}{\rm MQ} = 1.\end{eqnarray}
Thus
\begin{eqnarray}\label{two_a}
\int_{B_\epsilon(x)} \itfmq &=& \int_{({\rm Id}, tf)B_\epsilon(x)}
{\rm MQ}_{\Delta}
= \pm \int_{p_2  ({\rm Id}, tf)B_\epsilon(x)} (p_2^{-1})^*
{\rm MQ}_\Delta \\
&=& \pm\int_{(tf)(B_\epsilon(x))} {\rm MQ}_{\tmup} = \pm
\int_{\exp_x^{-1} (tf)(B_\epsilon(x))} \exp_x^*
{\rm MQ}_{\tmup}.\nonumber\end{eqnarray}
The last step is valid since for each $t$, 
$\exp_x^{-1}$ is well defined except on $tf(B_\epsilon(x)) \cap
\cutx$, which has measure zero.
By (\ref{one_a}), (\ref{two_a}), 
\begin{equation} \label{two_b} 
L(f) = \chi(M) \pm \lim_{\epsilon\to 0}\lim_{t\to 0} \int_{\exp_x^{-1}
(tf)(B_\epsilon(x))} \exp_x^* \MQ_{\tmup}.\end{equation}

We  modify this equation to handle the non-uniformity of the integral.
Since $\Gamma$ is transverse to $M\times \{f(x)\}$, for a fixed
$\epsilon^\prime <\epsilon
$ there exists $\delta = \delta(\epsilon^\prime)$ such that any
$\delta$ perturbation of $\Gamma$ in the $C^1$ topology is still a
graph over a set $U_{\epsilon^\prime}\subset \{x\}\times M$
containing $(x,f(x)).$  Also,
for any sequence $\epsilon_n \to 0$, there exists a sequence $t_n \to
0$ such that the graph $\Gamma_{t_n}$ of $t_nf$ is a $\delta_n =
\delta(\epsilon_n)$
perturbation of $\Gamma$.  Thus there is a set $U_n\subset
(t_nf)(B_{\epsilon_n}(x))$ such
that $\Gamma_{t_n}$ is a graph over $U_n.$  Set 
$W_n = (t_nf)^{-1}(U_n) \cap B_{\epsilon_n}(x)$.
Then
\begin{eqnarray} \label{inter}
\lim_{\epsilon\to 0} \lim_{t\to 0}
\int_{B_\epsilon(x)}({\rm Id}, tf)^* \MQ_{\Delta} &=& \lim_{n\to\infty}
\lim_{t\to 0}\biggl[ \int_{B_{\epsilon_n}(x) \setminus W_n} \itfmq
\nonumber\\
&&\qquad + 
\int_{ W_n} \itfmq \biggr]\nonumber\\
&=& \lim_{n\to\infty}\int_{B_{\epsilon_n}(x)\setminus
W_n}
\euler \\
&& \qquad + \lim_{n\to \infty}\lim_{t \to 0} 
\int_{ W_n} \itfmq\nonumber\\
&=& \lim_{n\to \infty}\lim_{t \to 0} 
\int_{ W_n} \itfmq.\nonumber\end{eqnarray}

The next technical lemma 
replaces $tf$ by a family of maps deforming $f(y)$
towards $x$ rather than towards $y$, for $x,y$ close.
\begin{lemma}\label{lemmayah}
 For $\mu >0$, there exists a neighborhood $U = U_\mu$ of
$x$ such that for all $y_0\in U$, there exists a unique minimal
geodesic 
$\gamma_{f(y_0),x}$ from
$f(y_0)$ to $x$ which is $\mu$ close in the
$C^1$ topology to the unique minimal geodesic 
$\gamma_{f(y_0),y_0}$ from $f(y_0)$ to $y_0$.  
\end{lemma}

\noindent{\sc Proof:}  The lemma is obvious unless $f(y_0)\in \cutx.$
In general, fix $y_0$ close to $x$ and let $y$ denote a point on the
minimal geodesic from $y_0$ to $x$. Since $f(y_0)\not\in {\cal
C}_{y_0}$, we have $y_0\not\in {\cal C}_{f(y_0)}$, and in particular
$y_0$ is not in the conjugate locus of $f(y_0)$.  Thus the exponential
map $\exp_{f(y_0)}:T_{f(y_0)}M\to M$ surjects onto some neighborhood
of $y_0$.  For $y$ close to $y_0$, there is a unique minimal
geodesic 
$\gamma_{f(y_0),y}$ from $f(y_0)$ to $y$, and the family of such
geodesics is $C^1$ close.   Now take
a curve $\gamma_\epsilon$ which is a smoothed approximation to
$\gamma_{f(y_0),y}$ followed by the minimal geodesic from $y$ to $x$
such that the length of $\gamma_\epsilon$ satisfies
$\ell(\gamma_\epsilon) \leq d(f(y_0),y) + d(y,x) + \epsilon.$
Parametrizing all curves by arclength, we see that for $y_0$ close
enough to $x$, the new family of curves is still $C^1$ close.  By the
Ascoli theorem, a subsequence of this family converges in $C^0$ as
$y\to x$ and as $\epsilon \to 0$
to a curve $\gamma_{f(y_0),x}$ from $f(y_0)$ to $x$ of length
$d(f(y_0),x)$--i.e. $\gamma_{f(y_0),x}$ is a minimal geodesic from
$f(y_0)$ to $x$.  Since
$\gamma_{f(y_0),x}$ is smooth and since the tangent vectors
$\exp_{f(y_0)}^{-1}y$ lie on the unit sphere in $T_{f(y_0)}M$, it
follows easily that a subsequence of  $\exp_{f(y_0)}^{-1}y$ converges
to a vector $v$ with $\exp_{f(y_0)}(sv) = \gamma_{f(y_0),x}.$  By the
smooth dependence of geodesics on initial conditions,  the minimal
geodesic $\gamma_{f(y_0),x}$ is $C^1$ close to the minimal geodesic 
from $f(y_0)$ to $y$, and hence $C^1$ close to $\gamma_{f(y_0),y_0}.$
 
This shows that along any radial geodesic $r$ centered at $x$, there
exists a distance $\delta = \delta(r)$ such that if $y$ is a point on
$r$ with $d(x,y) <\delta$, then there is a minimal geodesic from
$f(y)$ to $x$ which is $\mu$ close to the minimal geodesic from $f(y)$
to $y$ in the $C^1$ topolgy.  A similar argument shows that we may
take $\delta$ to be a continuous function of the radial direction.
\hfill$\Box$\break

\bigskip

We now fix $n$ large enough so that the lemma applies to all
$y\in W_n.$ Define a family of maps $g_t:W_n\to M$, $t\in (0,1]$,
which are approximations to $tf$ as follows.  For $x\in \cutf,$ set
$g_t(x) = (tf)(x) = x.$  For $y\not\in \cutf$ and $v_y =
\exp_{f(y)}^{-1} y,$  define $\alpha_t$ by $(tf)(y) =
\exp_{f(y)}(\alpha_t v_y)$.  Now set
$g_t(y) = \exp_ {f(y)}(\alpha_t v_x)$, where $\exp_{f(y)} (sv_x)$ is
the minimal geodesic from $f(y)$ to $x$ just constructed.  By the
smooth dependence of geodesics on parameters, we see that for $n$
large enough, $tf$ is arbitrarily $C^1$ close to $g_t$ for all $y\in
W_n$ and for all $t\in (0,1].$  This implies
\begin{eqnarray} \label{fg}
L(f) - \chi(M) &=&\lim_{n\to\infty}\lim_{t\to 0} \int_{W_n}
\itfmq\nonumber\\
&=&
\lim_{n\to\infty}\lim_{t\to 0} \int_{W_n} ({\rm
Id},g_t)^*{\MQ_\Delta}\nonumber\\
&&\qquad +
\lim_{n\to\infty}\lim_{t\to 0} \int_{W_n} [({\rm Id},tf)^* -({\rm
Id},g_t)^*] \MQ_\Delta\nonumber\\
&=& \lim_{n\to\infty}\lim_{t\to 0} \int_{W_n} ({\rm
Id},g_t)^*{\MQ_\Delta}\\
&&\qquad +
\lim_{n\to\infty}\int_{W_n} \lim_{t\to 0} [({\rm Id},tf)^* -({\rm
Id},g_t)^*] \MQ_\Delta\nonumber\\
&=& \lim_{n\to\infty}\lim_{t\to 0} \int_{W_n} ({\rm
Id},g_t)^*\MQ_\Delta.\nonumber \end{eqnarray}
As in (\ref{two_a}) we have
\begin{eqnarray}\label{up}
\int_{W_n} ({\rm
Id},g_t)^*\MQ_\Delta &=& \pm\int_{p_2({\rm Id},g_t)W_n} (p_2^{-1})^*
\MQ_\Delta\nonumber\\
&=& \pm\int_{\beta\exp_\cdot^{-1}g_t(W_n)}\MQ\\
&\approx& \pm\int_{\beta\exp_x^{-1}g_t(W_n)}\MQ,\nonumber
\end{eqnarray}
where $\exp_{\cdot}$ denotes the exponential map from
$TM^\uparrow|_{W_n}$
to $M\times M$. Since
$\exp_{q}^{-1}$ is $C^1$ close to $\exp_x^{-1}$ for $q$
 close to $x$, the error in the last line goes to zero as
$n\to\infty, t\to 0.$  (We use parallel translation to
compare maps with different ranges.)  Since
$\beta\exp_x^{-1}  g_tW_n\subset T_xM^\uparrow$, 
\begin{equation}\label{four_c}
 \biggl|\lim_{n\to\infty}\lim_{t\to 0} \int_{\beta\exp_{x}^{-1}
g_tW_n}\MQ\biggr| \leq \lim_{n\to\infty}\lim_{t\to 0}
\int_{T_xM^\uparrow} \biggl| \MQ\biggr| = 1,\end{equation}
as $\MQ$ is a positive multiple of the volume form in
$T_xM^\uparrow.$
Thus by  (\ref{fg})--(\ref{four_c}), we get
\begin{equation}\label{five} |L(f) -\chi(M) | \leq 1.\end{equation}

Summing over the finite number of points in $\cutf$ gives the main
theorem.
\begin{theorem} \label{isolated}
Assume that f is transverse to the cut locus.  Then
\[
|L(f) - \chi(M) | \leq |\cutf|.\]\end{theorem}

\bigskip

\noindent {\bf Remarks:}
1) The theorem is trivially sharp by setting $f = {\rm Id}.$  The
result
is also sharp for $f:z\mapsto z^n$ on $S^1$, and $f$  is transverse to the
cut locus.  
The two point suspension
of $f$ to $S^2$ is transverse to the cut locus and gives equality
in Theorem \ref{isolated},
and iterating this procedure gives sharp maps in all dimensions.

2) The inequality in Theorem \ref{isolated} can be refined to an
equality for maps on $S^n$ with the standard metric, since $\cutx =
\{-x\}$ easily implies $\lim_{n\to\infty}\lim_{t\to
0}\beta \exp_x^{-1} ({\rm Id}, g_t)W_n = T_xM^\uparrow$.  Thus as in
(\ref{yah}) the left hand side of (\ref{four_c}) is $\pm1$.  
If we denote $+1\
(-1)$ by sqn${}_x$ for $x\in \cutf$ if $p_2$ is orientation preserving
(reversing) on $\Gamma$ at $(x,f(x))$, then we have
\[ L(f) - \chi(M) = \sum_{x\in \cutf} {\rm sgn}_x.\]

\bigskip

We now show that the existence of a function $f$ transverse to the cut
locus and with $L(f) \neq \chi(M)$  imposes strong restrictions on the
metric. 
Recall that $\exp_x^{-1}g_t(y)$ lies on the radial line joining
$\exp_x^{-1} y$ to $0$ in $T_xM.$
Looking back at (\ref{four_c}), we have
\[\biggl|\lim_{n\to\infty}\lim_{t\to 0} \int_{\beta\exp_{x}^{-1}({\rm
Id},g_t)W_n}\MQ\biggr| = 1\]
iff $\exp_x^{-1}({\rm Id}, f)B_\epsilon(x)$ contains an interior
collar of the cut
locus in $T_xM$, as only in this case will
$\lim_{n\to\infty}\lim_{t\to
0}\beta\exp_x^{-1} ({\rm Id}, g_t)W_n = T_xM^\uparrow$. Letting
$\epsilon$
shrink, we see that this collar condition occurs only if the cut locus
of $x$ in $M$ is contained in an arbitrary neighborhood of
$f(x)$--i.e. the cut locus of $x$ in $M$ is precisely $f(x)$.  Thus
$M$ is homeomorphic to the one point compactification of $\exp_x^{-1}
(M\setminus \{f(x)\})$, so $M \approx S^n.$  

\begin{corollary}  (i)  Let $f:M\to M$ be a smooth map which is
transverse to the cut locus.  If $M\not\approx S^n$ and $\cutf \neq
\emptyset$,  then
\[ |L(f) - \chi(M) | < |\cutf|.\]

(ii) 
Let $f:M\to M$ be a smooth map which is Lefschetz and 
transverse to the cut locus.  Let ${\rm Fix}(f)$ be the fixed point
set of $f$.  
Then 
\[|\chi(M)| \leq |{\rm Fix}(f)| + |\cutf|,\]
with strict inequality if  $M\not\approx S^n$ and $\cutf \neq
\emptyset$.  \end{corollary}

\noindent (ii) follows from $|\chi(M)| \leq |L(f) - \chi(M)| + |L(f)|$
and the
Lefschetz fixed point theorem.  As an application of (i), let $f$ be
transverse to the cut locus and have
$L(f) \neq \chi(M)$.  Then if
$M\not\approx S^n$, either $\cutf 
=\emptyset$ or $|\cutf| \geq 2.$  This fails on $S^n$ for the
suspension of $z\mapsto z^2.$  

\bigskip

As in (\ref{four_c}),
set 
\begin{equation} \alpha_{x} = ({\rm sgn}\
p_2)\lim_{n\to\infty}\lim_{t\to 0} \int_{\beta\exp_x^{-1} g_t W_n}
 {\rm MQ}.\end{equation}
Then $\alpha_x \in [-1,1]$, and
$\sum_i \alpha_{x_i}\in L(f) - \chi(M) \in {\bf Z}$
by (\ref{fg}), (\ref{up}).  In general,
we expect $\alpha_x$ to vanish.  For let $A_\epsilon$ be the
radial projection of $\exp_x^{-1} f(B_\epsilon(x))$ onto $\cutx^\prime
 = \exp_x^{-1}\cutx$ in
$T_xM.$  Since the Mathai-Quillen form is a radially symmetric
fractional, positive
multiple of the volume form in each fiber, we have
$$|\alpha_x| \leq \frac{\lim_{\epsilon\to 0} |A_\epsilon|}{|\cutx^\prime|},$$
where $|A_\epsilon|,|\cutx^\prime|,$ denote the $(n-1)$-dimensional measure
of $A_\epsilon,\cutx^\prime.$  (Alternatively, we can project $\exp_x^{-1}
f(B_\epsilon(x))$ onto 
a small sphere centered at $0\in T_xM$ and compare its measure to the
measure of the sphere.)   Thus if $f$ is transverse to the cut locus,
and $L(f) \neq \chi(M)$, there exists $x\in \cutf$ such that
\begin{equation}\label{yahtwo}
\frac{\lim_{\epsilon\to 0} |A_\epsilon|}{|\cutx^\prime|} \neq 0.\end{equation}

By an application of Lemma \ref{lemmayah}, for $\epsilon$ small, for
all $y\in B_\epsilon(x)$, every minimal geodesic from $f(y)$ to $x$ is
$C^1$ close to a minimal geodesic from $f(x)$ to $x$.  Thus the unique
minimal geodesic $\gamma$ from $f(y)$ to $f(x)$ lifts under $\exp_x^{-1}$ to a
curve of length at most a constant times the length of $\gamma$, where
the constant depends on sectional curvature bounds for $M$.  Therefore
there exist $y'\in \exp_x^{-1}(f(y)),
\ z\in \exp_x^{-1}(f(x))$ such that $d( y',z) \to 0$ uniformly in $y$
as $\epsilon\to 0.$
It follows that (\ref{yahtwo}) can occur 
only if $\lim_{\delta\to 0}|B_\delta(\exp_x^{-1}
f(x))| \neq 0$, where $B_\delta$ denotes the delta ball in $\cutx^\prime$.
This
implies that $\exp_x^{-1} f(x)$ has positive $(n-1)$ dimensional
Hausdorff measure.

\begin{proposition}  \label{propend} Assume that
$f$ is transverse to the cut locus and $L(f) \neq \chi(M)$.  Then
there exists $x\in \cutf$ such that $\exp_x^{-1} f(x)$ has positive
$(n-1)$ dimensional Hausdorff measure in $\cutx^\prime.$ \end{proposition}

We call a metric on $M$ {\it somewhere (nowhere) sphere-like} if there
exist (do not exist)
$x,y\in M$ such that $\exp_x^{-1}y$ has positive 
$(n-1)$ dimensional Hausdorff measure in $\cutx^\prime.$  
One would expect a typical metric to be nowhere sphere-like, but
it is easy to construct
a somewhere sphere-like metric on any $M$, by
considering $M$ as the connect sum $M\# S^n$. 

\begin{theorem}\label{yahthree} (i) A metric of non-positive curvature
is nowhere sphere-like.

(ii)
Let $f:M\to M$ be a diffeomorphism
with $L(f) \neq \chi(M).$  If $M$ is nowhere sphere-like, then
$|\cutf| = \infty.$  \end{theorem}

\noindent{\sc Proof:}  (i) The map $\exp_x$ is a covering map for
metrics of non-positive curvature, so the inverse image of any $y\in M$
is discrete in $T_xM.$  

(ii)
By the proposition, $f$ is not transverse to
the cut locus.  Since condition (ii) for this 
transversality is satisfied, (i) must fail.  Thus $|\cutf| = \infty.$
\hfill$\Box$\break

\bigskip
\noindent {\bf Examples:}
(i) The flat torus and a constant negative curvature
surface are nowhere sphere-like, so no self-map of these spaces
satisfies the hypothesis of Proposition \ref{propend}.
For example,
for $x = (\theta,\psi) \in T^2
= [0,2\pi]\times [0,2\pi]$, the cut locus of $x$ in $T^2$ is
$(\{\theta\pm\pi\}\times [0,2\pi]) \cup ([0,2\pi] \times \{\psi\pm
\pi\}).$ 
For $(n,m)\in {\bf Z}^2$, 
$f(\theta,\psi) = (n\theta, m\psi)$ is a local diffeomorphism with
$L(f) = 2 - n - m.$  Theorem \ref{yahthree} applies to local
diffeomorphisms, so for $n+m\neq 2$, we conclude that
$|\cutf|=\infty.$  In fact, it is easy to check that $\cutf$ is
uncountable for $(n,m) \neq (1,1).$

(ii) Let $\Sigma^g$ be a genus $g>1$ surface symmetric about a plane
passing through the $g$ holes. For
$f:\Sigma^g\to\Sigma^g$ the diffeomorphism given by reflection
through this plane, $L(f) = 0.$
For example, $\Sigma^g$ can be the hyperelliptic
curve $y^2 = \prod_{i=1}^{2g+1} (x-a_i)$ with $a_i$ real and distinct,
with $f$ the involution $(x,y) \mapsto (x,-y)$.
  For any
metric on $\Sigma^g$ we have either (i) $f$ is not transverse to the
cut locus and so $|\cutf| = \infty$, or (ii) the metric is somewhere
sphere-like and $|\cutf| >2g-2.$  Thus for any metric, $|\cutf| >2g-2,$
and for most metrics $\cutf$ is infinite.  

(iii) Let $f:S^2\to S^2$ be a holomorphic map of degree
$n$.  Then $f$ is a branched covering and so $\Gamma$ is transverse to
$M\times \{f(x)\}$ except at the branch points $B$.  As above,
for any metric on $S^2$ with $f^{-1}(B)\cap\cutf = \emptyset$,
we have $|\cutf| \geq |L(f)-\chi(S^2)| = n-1.$

\bigskip
In \S3, the topological limit as $t\to\infty$ was shown
to have a
geometric refinement.  We do not have a
topological interpretation for the geometric $t\to 0$ limit;
presumably, the
singular part of the current of  \S4.1 represents 
a cohomology class in some  theory.

\appendix

\section{\large\rm\bf  Hodge theoretic techniques }

As mentioned in \S2.3, the upper bound for the Lefschetz number of a
flat manifold can be extended to arbitrary metrics.
Using sectional curvature bounds to control the Jacobi fields
and the curvature tensor, one can extract 
an upper bound from the integral formula Theorem \ref{localtheorem} in
terms of the sectional curvature.  In contrast, there is an easier
Hodge theory argument which
 constructs a better upper bound in terms of
Ricci curvature.

Let ${\cal N}= {\cal N}(n,C,D,V)$ be the class of Riemannian
$n$-manifolds $(M,g)$ with Ricci curvature Ric $\ge C$, diam$(M)\le
D$ and vol$(M)\ge V$.
\begin{proposition} There exist constants $C=C(k,n)$ and $D=D({\cal
N})$
such that for all $(M,g)\in {\cal N}$,
$$|L(f)|\le 1+D \sum_{k=1}^n C(k,n) \left( {n\atop
k}\right)\beta_k \cdot\sup_{x\in M} | df_x|_\infty^k,
$$
where $\beta_k$ is the $k^{th}$ Betti number of $M$.
\end{proposition}

Before the proof, we compare two norms for differential forms.
For $\alpha  \in\Lambda^k T_x^\ast M$, we have the
$L^2$ (Hodge) norm $|\alpha |_{2}^2 =\ast
(\alpha\land\ast\alpha)$ and the $L^\infty$ norm
$$|\alpha |_\infty =\sup_{v\in (T_xM)^{\otimes k} \setminus \{ 0\}}
\frac{|\alpha(v)|}{|v|},
$$
where $v= v_1\otimes \ldots\otimes v_k$ has norm $|v| = \prod |v_i|.$
Here we consider $\alpha$ as a linear functional on $(T_xM)^{\otimes
k}.$
Of
course, there exists $C=C(g)$ such that $C^{-1}|\alpha |_\infty\le
|\alpha|_{2} \le C|\alpha|_\infty$, but we want this constant to
depend only on $k,n$.
\begin{lemma} There exists a constant $C=C(k,n)$ such that
$$\left( {n\atop k}\right)^{-1/2} |\alpha|_{2}\le
|\alpha|_\infty \le C(k,n)|\alpha|_{2}.
$$
\end{lemma}

\noindent{\sc Proof:} Let $\{ \theta^i\}$ be an orthonormal basis of
$T_x^\ast
M$ with dual basis $\{ X_i\}$ of $T_xM$. For $\alpha
=\alpha_I\theta^I$, we have
$$|\alpha|_\infty \ge \frac{|(\alpha_I \theta^I)(X_{i_1}\otimes\ldots
\otimes X_{i_k})|}{| X_{i_1}\otimes\ldots\otimes X_{i_k}|}
=|\alpha_{I_0}|,$$
where $I_0=(i_1,\ldots, i_k)$. Thus
$$|\alpha|_\infty \ge \sup_I |\alpha_I|\ge
\left( {n\atop k}\right)^{-1/2}\left( \sum_I
|\alpha_I|^2\right)^{1/2} =\left( {n\atop
k}\right)^{-1/2} |\alpha|_{2}.
$$

For the other estimate
$$|\alpha|_\infty^2 \le \sup_{v=v_1\otimes\ldots\otimes
v_k\not= 0}\frac{\sum_I
|\alpha_I|^2 |\theta^I (v_1\otimes\ldots\otimes
v_k)|^2}{|v_1\otimes\ldots\otimes v_k|^2}.
$$
For fixed $I_0=(i_1,\ldots,i_k)$ and $v_1=a_1^{j_1}X_{j_1},\ldots,
v_k=a_k^{j_k} X_{j_k}$, we have
$$|\theta^{I_0} (v_1\otimes\ldots\otimes
v_k)|\le \sum_{\scriptstyle{{{j_1,\ldots,
j_k}\atop {{\{ j_1,\ldots,j_k\} =I_0}}}}} |a_1^{j_1}\cdot\ldots\cdot
a_k^{j_k}|.
$$
Thus 
\begin{eqnarray*}
|\alpha|_\infty^2 &\le& \sup_{v\not= 0}
\frac{\sum_{I_0}|\alpha_{I_0}|^2 \sum_{\{
j_1,\ldots,j_k\}=I_0}|a_1^{j_1}\cdot \ldots\cdot a_k^{j_k}|^2\cdot
k!}{|v_1\otimes\ldots\otimes v_k|^2} \\
&=& \sup_{v\not= 0}\frac{\sum_{I_0}|\alpha_{I_0}|^2\sum_{\{
j_1,\ldots,j_k\} =I_0}|a_1^{j_1}\cdot \ldots\cdot a_k^{j_k}|^2\cdot
k!}
{\prod_{q=1}^k (\sum_{l_q}(a_q^{l_q})^2)} \\
&=&\sup_{v\not= 0} \sum_{I_0}|\alpha_{I_0}|^2\left[ \frac{k!}{({n\atop
k})}\frac{\sum_{\{
j_1,\ldots,j_k\} =I_0}|a_1^{j_1}\cdot \ldots \cdot a_k^{j_k}|^2}
{\prod_{q=1}^k (\sum_{l_q}(a_q^{l_q})^2)}\right].
\end{eqnarray*}
For fixed $I_0$, the term inside the square bracket is a scale
invariant function on $\bbR^{nk}=\{ (a_i^j): i=1,\ldots,k,\
j=1,\ldots,
n\}$ and so is bounded above by $\bbC^\prime (k,n)$ independent of
$I_0$. Thus
$$|\alpha |_\infty^2\le \frac{k!}{({n\atop k})}C^\prime (k,n)\sum_I
|\alpha_I|^2 =(C(k,n))^2|\alpha |_{2}^2.$$
\hfill$\Box$\break

\bigskip

\noindent{\sc Proof of the Proposition:} Let $\{ \omega_k^i\}$  be an
$L^2$-orthonormal basis of harmonic
$k$-forms. The trace of $f^\ast : H^k(M;\bbR )\rightarrow H^k(M;\bbR)$
is
$\sum_i \langle f^\ast
\omega_k^i,\omega_k^i\rangle,$
so
\begin{eqnarray}\label{one}
|L(f)|&\le& \sum_{k,i}\biggl| \langle f^\ast \omega_k^i,
\omega_k^i\rangle\biggl| \le
\sum_{k,i}\Vert f^\ast \omega_k^i\Vert\nonumber\\
&=& \sum_{k,i} \left[ \int_M |(f^\ast \omega_k^i)_x|_2^2\mbox{\rm
dvol}(x)\right]^{1/2},
\end{eqnarray}
by Cauchy-Schwarz. 
Here $\Vert\alpha\Vert^2 = \int_M\alpha\wedge\ast\alpha$ is the global
$L^2$ norm. 
When $k=0$, we have $\Vert f^\ast \omega_0^1\Vert
=\Vert \omega_0^1\Vert =1$.

By (\ref{one}) and the lemma, we have
\begin{equation}\label{two}|L(f)| \le 1+\sum_{k=1}^n \sum_i\left(
{n\atop k}\right){\rm vol}^{1/2}
(M)\sup_{x\in M}|(f^\ast \omega_k^i)_x|_\infty.\end{equation}
Now
$$|(f^\ast \omega)_x|_{\infty} =\sup_{v\not= 0}
\frac{|(f^\ast \omega)_x
(v_1\otimes\ldots\otimes v_k)|}{|v_1\otimes\ldots\otimes 
v_k|}=\sup_{v\not=
0} \frac{|\omega_{f(x)}(f_\ast v_1\otimes\ldots\otimes f_\ast
v_k)|}{|v_1\otimes \ldots\otimes v_k|},
$$
where $f_* = df.$
Since the last term vanishes if $f_\ast v_i=0$ for some $i$, we assume
$f_\ast v_i\not= 0$. Then
\begin{eqnarray*}
|(f^\ast \omega)_x|_\infty &=& \sup_{v\not=
0}\frac{|\omega_{f(x)}(f_\ast
v_1\otimes\ldots\otimes f_\ast v_k)|}{|f_\ast v_i\otimes
\ldots\otimes f_\ast v_k|}\cdot
\frac{|f_\ast v_1\otimes\ldots\otimes f_\ast
v_k|}{|v_1\otimes\ldots\otimes v_k|} \\
&\le &|\omega_{f(x)}|_\infty\cdot \sup_{v\not= 0}
\frac{\prod_i|f_\ast v_i|}
{\prod_i|v_i|} \\
&\le &|\omega_{f(x)}|_\infty \cdot\sup_{v\not= 0} \frac{\prod_i |
df_x|_\infty |v_i|}{\prod_i|v_i|} \\
&\le &|\omega_{f(x)}|_\infty | df_x|_\infty^k.
\end{eqnarray*}
By (\ref{two}) and the lemma, we get
\begin{eqnarray*}|L(f)| &\le& 1+\sum_{k=1}^n\left( {n\atop
k}\right)\mbox{\rm
vol}^{1/2}(M)\sum_i \sup_{x\in M}| df_x|_\infty^k \cdot
|(\omega_k^i)_{f(x)}|_\infty.\\
&\leq& 1+\sum_{k=1}^n\left( {n\atop k}\right)\mbox{\rm
vol}^{1/2}(M)\sum_i \sup_{x\in M}| df_x|_\infty^k \cdot
C(k,n)|(\omega_k^i)_{f(x)}|_2.
\end{eqnarray*}
By \cite{C}, \cite{L}, 
there is an explicit constant $D_1({\cal N})$ such that for all
$x\in M$,
$$|(\omega_k^i)_x|_2\le D_1({\cal N})\Vert
\omega_k^i\Vert =D_1({\cal N}).
$$
Thus
$$|L(f)|\le 1+\sum_{k=1}^n \left( {n\atop
k}\right)\beta_k \cdot {\rm vol}^{1/2} (M)\cdot D_1({\cal N})C(k,n)
\sup_{x\in M}| df_x|_\infty^k.
$$
Finally, vol$(M)$ is bounded above on ${\cal N}$ by standard
comparison
theorems.\hfill$\Box$\break

\bibliographystyle{amsplain}
\bibliography{art5}

\end{document}